\documentclass[11pt]{amsart}
\usepackage[english]{babel}
\usepackage{amsthm} 
\usepackage{amssymb}
\usepackage{amsmath}
\usepackage{hyperref}
\usepackage{amsfonts}
\usepackage{verbatim}
\usepackage{color}
\usepackage{stmaryrd}
\usepackage{fancyhdr}
\usepackage{mathrsfs}
\usepackage{graphicx}
\usepackage{pdfsync}

\usepackage{esint} 

\newtheorem{theorem}{Theorem}[section]
\newtheorem{corollary}[theorem]{Corollary}
\newtheorem{definition}[theorem]{Definition}
\newtheorem{lemma}[theorem]{Lemma}
\newtheorem{proposition}[theorem]{Proposition}

\numberwithin{equation}{section}

\newcommand{\N}{\mathbb{N}}
\newcommand{\nn}{\mathbb{N}}
\newcommand{\R}{\mathbb{R}}
\newcommand{\rr}{\mathbb{R}}
\newcommand{\cc}{\mathbb{C}}
\newcommand{\CC}{\mathbb{C}}
\newcommand{\eps}{\varepsilon}

\newcommand{\real}{\mathbb{R}}
\def\un{{\mathrm{1~\hspace{-1.4ex}l}}}

\def\un{{\mathrm{1~\hspace{-1.4ex}l}}}

\usepackage{varioref}
 
\def\Rom#1{\uppercase\expandafter{\romannumeral #1}}

\def\polhk#1{\setbox0=\hbox{#1}{\ooalign{\hidewidth \lower1.5ex\hbox{`}\hidewidth\crcr\unhbox0}}}

\author{Karine \textsc{Beauchard}, Karel \textsc{Pravda-Starov}}

\address{\noindent \textsc{Karine Beauchard, IRMAR, \'Ecole normale sup\'erieure de Rennes, UBL, CNRS, Campus de Ker Lann, 35170 Bruz, France}}
\email{karine.beauchard@ens-rennes.fr}
\address{\noindent \textsc{Karel Pravda-Starov, IRMAR, CNRS UMR 6625, Universit\'e de Rennes 1, Campus de Beaulieu, 263 avenue du G\'en\'eral Leclerc, CS 74205,
35042 Rennes cedex, France}}
\email{karel.pravda-starov@univ-rennes1.fr}


\thanks{The authors were partially supported by the ``Agence Nationale de la Recherche'' ANR EMAQS (Project: ANR-2011-BS01-017-01) and ANR NOSEVOL 
(Project: ANR 2011-BS01-019-01).}

\title{Null-controllability of hypoelliptic quadratic differential equations}
\date{}

\begin{document}

\begin{abstract}
We study the null-controllability of parabolic equations associated with a general class of hypoelliptic quadratic differential operators. Quadratic differential operators are operators defined in the Weyl quantization by complex-valued quadratic symbols. We consider in this work the class of accretive quadratic operators with zero singular spaces. These possibly degenerate non-selfadjoint differential operators are known to be hypoelliptic and to generate contraction semigroups which are smoothing in specific Gelfand-Shilov spaces for any positive time. Thanks to this regularizing effect, we prove by adapting the Lebeau-Robbiano method that parabolic equations associated with these operators are null-controllable in any positive time from control regions, for which null-controllability is classically known to hold in the case of the heat equation on the whole space. Some applications of this result are then given to the study of parabolic equations associated with hypoelliptic Ornstein-Uhlenbeck operators acting on weighted $L^2$ spaces with respect to invariant measures. By using the same strategy, we also establish the null-controllability in any positive time from the same control regions for parabolic equations associated with any hypoelliptic Ornstein-Uhlenbeck operator acting on the flat $L^2$ space extending in particular the known results for the heat equation or the Kolmogorov equation on the whole space.    
\end{abstract}

\maketitle

\tableofcontents

\section{Introduction}

\subsection{Null-controllability of degenerate parabolic equations}

We aim in this work at studying the null-controllability of parabolic equations controlled by a source term $u$ locally distributed on an open subset $\omega \subset \mathbb{R}^n$ of the whole space
\begin{equation}\label{syst_general}
\left\lbrace \begin{array}{ll}
(\partial_t + P)f(t,x)=u(t,x)\un_{\omega}(x)\,, \quad &  x \in \mathbb{R}^n, t>0, \\
f|_{t=0}=f_0 \in L^2(\rr^n),                                       &  
\end{array}\right.
\end{equation}
where $P=q^w(x,D_x)$ is an accretive quadratic operator. Quadratic operators are pseudodifferential operators defined in the Weyl quantization
\begin{equation}\label{3}
q^w(x,D_x) f(x) =\frac{1}{(2\pi)^n}\int_{\rr^{2n}}{e^{i(x-y) \cdot \xi}q\Big(\frac{x+y}{2},\xi\Big)f(y)dyd\xi}, 
\end{equation}
by symbols $q(x,\xi)$, with $(x,\xi) \in \rr^{n} \times \rr^n$, $n \geq 1$, which are complex-valued quadratic forms 
\begin{eqnarray*}
q : \rr_x^n \times \rr_{\xi}^n &\rightarrow& \cc\\
 (x,\xi) & \mapsto & q(x,\xi).
\end{eqnarray*}
These operators are non-selfadjoint differential operators in general, with simple and fully explicit expression since the Weyl quantization of the quadratic symbol $x^{\alpha} \xi^{\beta}$, with $(\alpha,\beta) \in \nn^{2n}$, $|\alpha+\beta| = 2$, is the differential operator
$$\frac{x^{\alpha}D_x^{\beta}+D_x^{\beta} x^{\alpha}}{2}, \quad D_x=i^{-1}\partial_x.$$
We study the null-controllability of the parabolic equations (\ref{syst_general}) associated with a general class of hypoelliptic quadratic differential operators:

\medskip

\begin{definition} [Null-controllability] Let $T>0$ and $\omega$ be an open subset of $\mathbb{R}^n$. 
Equation \eqref{syst_general} is said to be {\em null-controllable from the set $\omega$ in time} $T$ if, for any initial datum $f_0 \in L^{2}(\mathbb{R}^n)$, there exists $u \in L^2((0,T)\times\mathbb{R}^n)$, supported in $(0,T)\times\omega$, such that the mild solution of \eqref{syst_general} satisfies $f(T,\cdot)=0$.
\end{definition}

\medskip

By the Hilbert Uniqueness Method, see \cite[Theorem~2.44]{coron_book} or \cite{JLL_book}, the null-controllability of the equation \eqref{syst_general} is equivalent to the observability of the adjoint system 
\begin{equation} \label{adj_general}
\left\lbrace \begin{array}{ll}
(\partial_t + P^*)g(t,x)=0\,, \quad & x \in \mathbb{R}^n\,, \\
g|_{t=0}=g_0 \in L^2(\rr^n).
\end{array}\right.
\end{equation}
We recall that the notion of observability is defined as follows:

\medskip

\begin{definition} [Observability] Let $T>0$ and $\omega$ be an open subset of $\mathbb{R}^n$. 
Equation \eqref{adj_general} is said to be {\em observable in the set $\omega$ in time} $T$ if there exists a constant $C_T>0$ such that,
for any initial datum $g_0 \in L^{2}(\mathbb{R}^n)$, the mild solution of \eqref{adj_general} satisfies
\begin{equation}\label{eq:observability}
\int\limits_{\mathbb{R}^n} |g(T,x)|^{2} dx  \leq C_T \int\limits_{0}^{T} \Big(\int\limits_{\omega} |g(t,x)|^{2} dx\Big) dt\,.
\end{equation}
\end{definition}

\medskip

An important open problem at the core of current investigations is to understand to which extent the null-controllability (or observability) results known for uniformly parabolic equations still hold for degenerate parabolic equations of hypoelliptic type.

For equations posed on bounded domains, some progress have been made. In the case of the heat equation on a bounded domain $\Omega$ with Dirichlet boundary conditions on $\partial \Omega$, it is well-known that observability holds in arbitrary positive time $T>0$, with any non-empty open set $\omega$, see~\cite[Theorem~3.3]{Fattorini-Russel}, \cite{Fursikov-Imanuvilov-186} and \cite{Lebeau-Robbiano}. Degenerate parabolic equations exhibit a wider range of behaviours. Indeed, observability may hold true, or not, depending on the strength of the degeneracy. This feature is well understood for parabolic equations that degenerate on the domain boundary, see~\cite{Ala-Can-Fra,Can-Fra-Roc_2,Can-Fra-Roc_1, Cannarsa-V-M-ADE, Cannarsa-V-M-SIAM, Martinez-Vancost-JEE-2006} in the one-dimensional case, and~\cite{Cannarsa-V-M-CRAS} for the multi-dimensional one. Furthermore, a positive minimal time may be required to get observability, see the works~\cite{Grushin, MR3162108} in the case of the Grushin equation, \cite{KB_Can_Heisenberg} for the Heisenberg heat equation, and~\cite{MR3163490} for the Kolmogorov equation. This minimal time is actually related to localization properties of eigenfunctions. Finally, a geometric control condition may also be required for the observability inequality to hold~\cite{KB_Helffer}.

On the other hand, the understanding of the null-controllability (or observability) for degenerate parabolic equations of hypoelliptic type posed on the whole space is still at an earlier stage. For the heat equation on the whole space
\begin{equation}\label{heat_ws}
\left( \partial_t - \Delta_x \right)f(t,x)=u(t,x)\un_{\omega}(x)\,, \quad (t,x)\in(0,T)\times\mathbb{R}^n,
\end{equation}
no necessary and sufficient condition on the control region $\omega$ is known for null-controllability to hold in any positive time. The condition
$$\sup_{x \in \mathbb{R}^n}d(x,\omega)<+\infty,$$
is shown in~\cite[Theorem 1.11]{Miller_BSM2005} to be necessary for null-controllability to hold in any positive time. On the other hand, the following sufficient condition  
\begin{equation}\label{hyp_omega}
\exists \delta, r >0, \forall y \in \mathbb{R}^n\,, \exists y' \in \omega,\quad B(y',r) \subset \omega \text{ and } |y-y'|<\delta\,,
\end{equation}
is given in~\cite{Miller_unbounded} for null-controllability to hold from the open set $\omega \subset \rr^n$ in any positive time. 
The very same condition is shown in~\cite{YubiaoZhang} to be sufficient for the null-controllability of the Kolmogorov equation
\begin{equation}\label{koleq}
\left\lbrace \begin{array}{ll}
(\partial_t + v \cdot \nabla_x-\Delta_v)f(t,x)=u(t,x)\un_{\omega}(x)\,, \quad &  x \in \mathbb{R}^n,\\
f|_{t=0}=f_0 \in L^2(\rr^n),                                       &  
\end{array}\right.
\end{equation}
to hold in any positive time, see also~\cite{LeRousseau_Moyano} for control sets with Cartesian product structures $\omega=\omega_x \times \omega_v$. This result relies on a key spectral inequality proved in~\cite{LeRousseau_Moyano}.

As a first result in this work (Theorem~\ref{meta_thm}), we prove that condition (\ref{hyp_omega}) is actually sufficient for the null-controllability of all hypoelliptic Ornstein-Uhlenbeck equations
\begin{equation} \label{syst_LR}
\left\lbrace \begin{array}{ll}
\partial_t f(t,x) - \frac{1}{2}\textrm{Tr}[Q\nabla_x^2 f(t,x)] - \langle Bx, \nabla_x f(t,x)\rangle=u(t,x)\un_{\omega}(x)\,, \quad &  x \in \mathbb{R}^n,\\
f|_{t=0}=f_0 \in L^2(\rr^n),                                       &  
\end{array}\right.
\end{equation}
where $Q$ and $B$ are real $n \times n$-matrices satisfying the Kalman rank condition, with $Q$ symmetric positive semidefinite. This general result allows one to recover in particular the results of null-controllability
for the heat equation and the Kolmogorov equation.
Our proof relies on an adaptation of the Lebeau-Robbiano strategy. Compared to the classical Lebeau-Robbiano method, the new difficulty in the present analysis is that the above Ornstein-Uhlenbeck semigroups do not commute with the Fourier frequency cutoff projections. In order to address this problem, we need to adapt the Lebeau-Robbiano strategy by taking advantage of some key Gevrey smoothing properties of the Ornstein-Uhlenbeck equations.

In the second part of this work, we prove that the parabolic equations (\ref{syst_general}) associated with a general class of hypoelliptic quadratic operators are null-controllable from open sets satisfying condition (\ref{hyp_omega}) in any positive time. More specifically, our main result (Theorem~\ref{Main_result_control}) establishes that null-controllability holds for the parabolic equation (\ref{syst_general}), as soon as the Weyl symbol $q$ of the quadratic operator $q^w(x,D_x)$ has a non-negative real part $\textrm{Re }q \geq 0$ and a zero singular space $S=\{0\}$. The notion of singular space was introduced in~\cite{kps2} by Hitrik and the second author by pointing out the existence of a particular vector subspace in the phase space $S \subset \rr^{2n}$, which is intrinsically associated with a quadratic symbol~$q$. As pointed out in \cite{kps2,short,HPSVII,kps11,kps21,rodwahl,viola0}, the notion of singular space plays a basic role in the understanding of the spectral and hypoelliptic properties of the (possibly) non-elliptic quadratic operator $q^w(x,D_x)$, as well as the spectral and pseudospectral properties of certain classes of degenerate doubly characteristic pseudodifferential operators~\cite{kps3,kps4,viola1,viola2}. In particular, the work~\cite[Theorem~1.2.2]{kps2} gives a complete description for the spectrum of any non-elliptic quadratic operator $q^w(x,D_x)$ whose Weyl symbol $q$ has a non-negative real part $\textrm{Re }q \geq 0$, and satisfies a condition of partial ellipticity along its singular space~$S$,
\begin{equation}\label{sm2}
(x,\xi) \in S, \ q(x,\xi)=0 \Rightarrow (x,\xi)=0. 
\end{equation}
Under these assumptions, the spectrum of the quadratic operator $q^w(x,D_x)$ is shown to be composed of a countable number of eigenvalues with finite algebraic multiplicities. The structure of this spectrum is similar to the one known for elliptic quadratic operators~\cite{sjostrand}. This condition of partial ellipticity is generally weaker than the condition of ellipticity, $S \subsetneq \rr^{2n}$, and allows one to deal with more degenerate situations. An important class of quadratic operators satisfying condition (\ref{sm2}) are those with zero singular spaces $S=\{0\}$. In this case, the condition of partial ellipticity trivially holds.
More specifically, these quadratic operators have been shown in \cite[Theorem~1.2.1]{kps21} to be hypoelliptic and to enjoy global subelliptic estimates of the type
\begin{equation}\label{lol1}
\exists C>0, \forall u \in \mathscr{S}(\rr^n), \quad \|\langle(x,D_x)\rangle^{2(1-\delta)} u\|_{L^2} \leq C(\|q^w(x,D_x) u\|_{L^2}+\|u\|_{L^2}),
\end{equation}
where $\langle(x,D_x)\rangle^{2}=1+|x|^2+|D_x|^2$, with a sharp loss of derivatives $0 \leq \delta<1$ with respect to the elliptic case (case $\delta=0$), which can be explicitly derived from the structure of the singular space. Our proof of null-controllability for the parabolic equation (\ref{syst_general}) associated with a quadratic operator $P=q^w(x,D_x)$, whose Weyl symbol $q$ has a non-negative real part $\textrm{Re }q \geq 0$ and a zero singular space $S=\{0\}$, relies on a similar adaptation of the Lebeau-Robbiano strategy as the one devised for Ornstein-Uhlenbeck equations  (\ref{syst_LR}). Contrary to the case of Ornstein-Uhlenbeck equations, the analysis of this class of quadratic operators that are differential operators with variable coefficients, cannot rely on a sole frequency analysis on the Fourier side. We actually use some recent microlocal results on the Gelfand-Shilov regularizing properties of the semigroups generated by these quadratic operators together with a spectral inequality for Hermite functions (Proposition~\ref{thm:Spect_Ineq}). As for Ornstein-Uhlenbeck equations, the main difficulty is that the above semigroups do not necessarily commute with the projections onto Hermite functions. In order to address this problem, we need to adapt the Lebeau-Robbiano strategy by taking advantage of some key Gelfand-Shilov smoothing properties of these semigroups.

\subsection{Miscellaneous facts about quadratic differential operators}

Let $q^w(x,D_x)$ be a quadratic operator defined by the Weyl quantization (\ref{3}) of a complex-valued quadratic form $q$ on the phase space $\rr^{2n}$.
The maximal closed realization of the quadratic operator $q^w(x,D_x)$ on $L^2(\rr^n)$, that is, the operator equipped with the domain
\begin{equation}\label{dom1}
D(q^w)=\big\{g \in L^2(\rr^n) : \ q^w(x,D_x)g \in L^2(\rr^n)\big\},
\end{equation}
where $q^w(x,D_x)g$ is defined in the distribution sense, is known to coincide with the graph closure of its restriction to the Schwartz space~\cite{mehler} (pp.~425-426),
$$q^w(x,D_x) : \mathscr{S}(\rr^n) \rightarrow \mathscr{S}(\rr^n).$$
Classically, to any quadratic form defined on the phase space 
$$q : \rr_x^n \times \rr_{\xi}^n \rightarrow \mathbb{C},$$
is associated a matrix $F \in M_{2n}(\CC)$ called its Hamilton map, or its fundamental matrix, which is defined as the unique matrix satisfying the identity
\begin{equation}\label{10}
\forall  (x,\xi) \in \R^{2n},\forall (y,\eta) \in \R^{2n}, \quad q((x,\xi),(y,\eta))=\sigma((x,\xi),F(y,\eta)), 
\end{equation}
with $q(\cdot,\cdot)$ the polarized form associated with the quadratic form $q$, where $\sigma$ stands for the standard symplectic form
\begin{equation}\label{11}
\sigma((x,\xi),(y,\eta))=\langle \xi, y \rangle -\langle x, \eta\rangle=\sum_{j=1}^n(\xi_j y_j-x_j \eta_j),
\end{equation}
with $x=(x_1,...,x_n)$, $y=(y_1,....,y_n)$, $\xi=(\xi_1,...,\xi_n)$, $\eta=(\eta_1,...,\eta_n) \in \cc^n$. We observe from the definition that 
$$F=\frac{1}{2}\left(\begin{array}{cc}
\nabla_{\xi}\nabla_x q & \nabla_{\xi}^2q  \\
-\nabla_x^2q & -\nabla_{x}\nabla_{\xi} q 
\end{array} \right),$$
where the matrices $\nabla_x^2q=(a_{i,j})_{1 \leq i,j \leq n}$,  $\nabla_{\xi}^2q=(b_{i,j})_{1 \leq i,j \leq n}$, $\nabla_{\xi}\nabla_x q =(c_{i,j})_{1 \leq i,j \leq n}$,
$\nabla_{x}\nabla_{\xi} q=(d_{i,j})_{1 \leq i,j \leq n}$ are defined by the entries
$$a_{i,j}=\partial_{x_i,x_j}^2 q, \quad b_{i,j}=\partial_{\xi_i,\xi_j}^2q, \quad c_{i,j}=\partial_{\xi_i,x_j}^2q, \quad d_{i,j}=\partial_{x_i,\xi_j}^2q.$$
The notion of singular space introduced in~\cite{kps2} by Hitrik and the second author is defined as the following finite intersection of kernels
\begin{equation}\label{h1bis}
S=\Big( \bigcap_{j=0}^{2n-1}\textrm{Ker}
\big[\textrm{Re }F(\textrm{Im }F)^j \big]\Big)\cap \rr^{2n},
\end{equation}
where $\textrm{Re }F$ and $\textrm{Im }F$ stand respectively for the real and imaginary parts of the Hamilton map $F$ associated with the quadratic symbol $q$,
$$\textrm{Re }F=\frac{1}{2}(F+\overline{F}), \quad \textrm{Im }F=\frac{1}{2i}(F-\overline{F}).$$

When the quadratic symbol $q$ has a non-negative real part $\textrm{Re }q \geq 0$, the singular space can be defined in an equivalent way as the subspace in the phase space where all the Poisson brackets 
$$H_{\textrm{Im}q}^k \textrm{Re }q=\left(\frac{\partial \textrm{Im }q}{\partial\xi}\cdot \frac{\partial}{\partial x}-\frac{\partial \textrm{Im }q}{\partial x}\cdot \frac{\partial}{\partial \xi}\right)^k \textrm{Re } q, \quad k \geq 0,$$ 
are vanishing
$$S=\big\{X=(x,\xi) \in \rr^{2n} : \ (H_{\textrm{Im}q}^k \textrm{Re } q)(X)=0,\ k \geq 0\big\}.$$
This dynamical definition shows that the singular space corresponds exactly to the set of points $X \in \rr^{2n}$, where the real part of the symbol $\textrm{Re }q$ under the flow of the Hamilton vector field $H_{\textrm{Im}q}$ associated with its imaginary part
\begin{equation}\label{evg5}
t \mapsto \textrm{Re }q(e^{tH_{\textrm{Im}q}}X),
\end{equation}
vanishes to infinite order at $t=0$. This is also equivalent to the fact that the function \eqref{evg5} is identically zero on~$\rr$.

In this work, we study the class of quadratic operators whose Weyl symbols have non-negative real parts $\textrm{Re }q \geq 0$, and zero singular spaces $S=\{0\}$. 
According to the above description of the singular space, these quadratic operators are exactly those whose Weyl symbols have a non-negative real part $\textrm{Re }q \geq 0$, 
becoming positive definite 
\begin{equation}\label{evg41}
\forall T>0, \quad \langle \textrm{Re }q\rangle_T(X)=\frac{1}{2T}\int_{-T}^T{(\textrm{Re }q)(e^{tH_{\textrm{Im}q}}X)dt} \gg 0, 
\end{equation}
after averaging by the linear flow of the Hamilton vector field associated with its imaginary part.
These quadratic operators are also known~\cite[Theorem~1.2.1]{kps2} to generate contraction semigroups $(e^{-tq^w})_{t \geq 0}$ on $L^2(\rr^n)$, which are smoothing in the Schwartz space for any positive time
$$\forall t>0, \forall g \in L^2(\rr^n), \quad e^{-t q^w}g \in \mathscr{S}(\rr^n).$$
In all the following, the terminology semigroup refers to strongly continuous one parameter semigroup. 
In the recent work~\cite[Theorem~1.2]{HPSVII}, these regularizing properties were sharpened and 
these contraction semigroups were shown to be actually smoothing for any positive time in the Gelfand-Shilov space $S_{1/2}^{1/2}(\rr^n)$: $\exists C>0$, $\exists t_0 > 0$, $\forall g\in L^2(\real^n)$, $\forall \alpha, \beta \in \nn^n$, $\forall 0<t \leq t_0$,
\begin{equation}\label{eq1.10}
\|x^{\alpha}\partial_x^{\beta}(e^{-tq^w}g)\|_{L^{\infty}(\rr^n)} \leq \frac{C^{1+|\alpha|+|\beta|}}{t^{\frac{2k_0+1}{2}(|\alpha|+|\beta|+2n+s)}}(\alpha!)^{1/2}(\beta!)^{1/2}\|g\|_{L^2(\rr^n)},
\end{equation}
where $s$ is a fixed integer verifying $s > n/2$, and where $0 \leq k_0 \leq 2n-1$ is the smallest integer satisfying 
\begin{equation}\label{h1bis2}
\Big( \bigcap_{j=0}^{k_0}\textrm{Ker}\big[\textrm{Re }F(\textrm{Im }F)^j \big]\Big)\cap \rr^{2n}=\{0\}.
\end{equation}
As mentioned above, this Gelfand-Shilov regularizing property will be a key ingredient for deriving observability estimates in Section~\ref{sec:Proof}.

A first interesting example of an accretive quadratic operator with a zero singular space $S=\{0\}$ is given by the Kramers-Fokker-Planck operator acting on $L^2(\rr_{x,v}^2)$,
\begin{equation} \label{eq:KFP}
K=-\Delta_v+\frac{v^2}{4}+v\partial_x-\nabla_xV(x)\partial_v, \quad (x,v) \in \rr^{2},
\end{equation}
with a quadratic potential
$$V(x)=\frac{1}{2}ax^2, \quad a \in \rr^*.$$
Indeed, this operator writes as $K=q^w(x,v,D_x,D_v)$, where
$$q(x,v,\xi,\eta)=\eta^2+\frac{1}{4}v^2+i(v \xi-a x \eta),$$
is a non-elliptic complex-valued quadratic form with a non-negative real part, whose Hamilton map is given by 
$$F= \left( \begin{array}{cccc}
  0 & \frac{1}{2}i & 0 & 0 \\
  -\frac{1}{2}ai& 0 & 0 & 1 \\
0 & 0 & 0 &  \frac{1}{2}ai \\
0 & -\frac{1}{4} & -\frac{1}{2}i & 0 
  \end{array}
\right).$$
A simple algebraic computation shows that 
\begin{equation} \label{k0:KFP}
\textrm{Ker}(\textrm{Re }F) \cap \textrm{Ker}(\textrm{Re }F \ \textrm{Im }F) \cap \rr^{4}=\{0\}.
\end{equation}
The singular space of $q$ is therefore equal to zero $S=\{0\}$.
For the Kramers-Fokker-Planck operator, the integer $0 \leq k_0 \leq 2n-1$ defined in (\ref{h1bis2}) with here $n=2$, is equal to $1$.

According to  \cite[Theorem~1.2.1]{kps21}, this integer $0 \leq k_0 \leq 2n-1$ is directly related to the loss of derivatives $0 \leq \delta=2k_0/(2k_0+1)<1$ in the global subelliptic estimate (\ref{lol1}) satisfied by any quadratic operator whose Weyl symbol has a non-negative real part and a zero singular space. 
The following examples show that this integer can actually take any value in the set $\{0,...,2n-1\}$, when $n \geq 1$:

\medskip

\begin{itemize}
\item[-] Case $k_0=0$: Any quadratic symbol $q$ with $\textrm{Re }q \gg 0$ a positive definite real part
\item[-] Case $k_0=1$:
$$q(x,\xi)=\xi_2^2+x_2^2+i(x_2 \xi_1-x_1 \xi_2)+\sum_{j=3}^n(\xi_j^2+x_j^2)$$
\item[-] Case $k_0=2p$, with $1 \leq p \leq n-1$: 
$$q(x,\xi)=\xi_1^2+x_1^2+i(\xi_1^2+2x_2\xi_1+\xi_2^2+2 x_3\xi_2+....+\xi_p^2+2x_{p+1}\xi_p+\xi_{p+1}^2)+\sum_{j = p+2}^n(\xi_j^2+x_j^2)$$
\item[-] Case $k_0=2p+1$, with $1 \leq p \leq n-1$: 
$$q(x,\xi)=x_1^2+i(\xi_1^2+2x_2\xi_1+\xi_2^2+2 x_3\xi_2+....+\xi_p^2+2x_{p+1}\xi_p+\xi_{p+1}^2)+\sum_{j = p+2}^n(\xi_j^2+x_j^2)$$
\end{itemize}

\subsection{Statements of the main results}

\subsubsection{Null-controllability of hypoelliptic Ornstein-Uhlenbeck equations}\label{sec:LRsub}
We consider Ornstein-Uhlenbeck operators 
\begin{equation}\label{jen0}
P=\frac{1}{2}\sum_{i,j=1}^{n}q_{i,j}\partial_{x_i,x_j}^2+\sum_{i,j=1}^nb_{i,j}x_j\partial_{x_i}=\frac{1}{2}\textrm{Tr}(Q\nabla_x^2)+\langle Bx,\nabla_x\rangle, \quad x \in \rr^n,
\end{equation}
where $Q=(q_{i,j})_{1 \leq i,j \leq n}$ and $B=(b_{i,j})_{1 \leq i,j \leq n}$ are real $n \times n$-matrices, with $Q$ symmetric positive semidefinite. 
We denote $\langle A,B \rangle$ and $|A|^2$ the scalar operators
\begin{equation}\label{not}
\langle A,B\rangle=\sum_{j=1}^nA_jB_j, \quad |A|^2=\langle A,A\rangle=\sum_{j=1}^nA_j^2,
\end{equation}
when $A=(A_1,...,A_n)$ and $B=(B_1,...,B_n)$ are vector-valued operators. Notice that $\langle A,B\rangle\neq\langle B,A\rangle$ in general, since e.g., 
$\langle \nabla_x,Bx\rangle=\langle Bx,\nabla_x\rangle +\textrm{Tr}(B).$

We study degenerate hypoelliptic Ornstein-Uhlenbeck operators when the symmetric matrix $Q$ is possibly not positive definite. 
These degenerate operators have been studied in the recent works~\cite{lanco,farkas,lunardi1,lanco2,Lorenzi,Metafune_al2002,uhlenb}. 
We recall from these works that the assumption of hypoellipticity is classically characterized by the following equivalent assertions:

\medskip
 
\begin{itemize}
\item[$(i)$] The Ornstein-Uhlenbeck operator $P$ is hypoelliptic
\item[$(ii)$] The symmetric positive semidefinite matrices
\begin{equation}\label{pav0}
Q_t=\int_0^{t}e^{sB}Qe^{sB^T}ds,
\end{equation}
with $B^T$ the transpose matrix of $B$, are nonsingular for some (equivalently, for all) $t>0$, i.e. $\det Q_t>0$
\item[$(iii)$] The Kalman rank condition holds: 
\begin{equation}\label{kal1}
\textrm{Rank}[B|Q^{\frac{1}{2}}]=n,
\end{equation} 
where 
$$[B|Q^{\frac{1}{2}}]=[Q^{\frac{1}{2}},BQ^{\frac{1}{2}},\dots, B^{n-1}Q^{\frac{1}{2}}],$$ 
is the $n\times n^2$ matrix obtained by writing consecutively the columns of the matrices $B^jQ^{\frac{1}{2}}$, with $Q^{\frac{1}{2}}$ the symmetric positive semidefinite matrix given by the square root of $Q$
\item[$(iv)$] The H\"ormander condition holds:
$$\forall x \in \rr^n, \quad \textrm{Rank } \mathcal{L}(X_1,X_2,...,X_n,Y_0)(x)=n,$$
with 
$$Y_0=\langle Bx,\nabla_x\rangle, \quad X_i=\sum_{j=1}^nq_{i,j}\partial_{x_j}, \quad i=1,...,n,$$
where $ \mathcal{L}(X_1,X_2,...,X_n,Y_0)(x)$ denotes the Lie algebra generated by the vector fields $X_1$, $X_2$, ..., $X_n$, $Y_0$, at point $x \in \rr^n$ 
 \end{itemize}

\medskip

\noindent
When the Ornstein-Uhlenbeck operator is hypoelliptic, that is, when one (equivalently, all) of the above conditions holds, 
the associated Markov semigroup $(T(t))_{t \geq 0}$ has the following explicit representation due to Kolmogorov~\cite{kolmo}:
\begin{equation}\label{gen}
(e^{tP}f)(x)=\frac{1}{(2\pi)^{\frac{n}{2}}\sqrt{\det Q_t}}\int_{\rr^n}e^{-\frac{1}{2}\langle Q_t^{-1}y,y\rangle}f(e^{tB}x-y)dy, \quad t>0.
\end{equation}
The first result contained in this work establishes the null-controllability of any hypoelliptic Ornstein-Uhlenbeck equation from any open control region satisfying condition (\ref{hyp_omega}) in any positive time:

\medskip

\begin{theorem} \label{meta_thm}
Let $T>0$ and $\omega$ be an open subset of $\mathbb{R}^n$ satisfying (\ref{hyp_omega}). When the Kalman rank condition (\ref{kal1}) holds, the Ornstein-Uhlenbeck equation posed in the $L^2(\rr^n)$ space
\begin{equation} \label{syst_LR1}
\left\lbrace \begin{array}{ll}
\partial_t f(t,x) - \frac{1}{2}\emph{\textrm{Tr}}[Q\nabla_x^2 f(t,x)] - \langle Bx, \nabla_x f(t,x)\rangle=u(t,x)\un_{\omega}(x)\,, \quad &  x \in \mathbb{R}^n,\\
f|_{t=0}=f_0 \in L^2(\rr^n),                                       &  
\end{array}\right.
\end{equation}
is null-controllable from the set $\omega$ in time $T>0$.
\end{theorem}

\medskip

This theorem allows one in particular while taking $Q=2I_n$ and $B=0$, to recover the result of null-controllability of the heat equation (\ref{heat_ws}) proved in~\cite{Miller_unbounded}, and by taking 
$$Q=\left(\begin{array}{ll}
0 & 0 \\ 0 & 2I_d
\end{array}\right), \qquad
B=\left(\begin{array}{ll}
0 & -I_d \\ 0 & 0
\end{array}\right), \qquad n=2d,$$
to also recover the result of null-controllability of the Kolmogorov equation (\ref{koleq}) proved in~\cite{YubiaoZhang}.
The proof of Theorem~\ref{meta_thm} is given in Section~\ref{sec:LR}.

\subsubsection{Null-controllability and observability of parabolic equations associated with accretive quadratic operators with zero singular spaces}\label{sec:LRquad}

The main result contained in this article is the following:

\medskip

\begin{theorem}\label{Main_result_control}
Let $q : \rr_{x}^{n} \times \rr_{\xi}^n \rightarrow \cc$ be a complex-valued quadratic form with a non negative real part $\emph{\textrm{Re }}q \geq 0$, and a zero singular space $S=\{0\}$. If  $\omega$ is an open subset of $\mathbb{R}^n$ satisfying condition (\ref{hyp_omega}), then the parabolic equation 
$$\left\lbrace \begin{array}{ll}
\partial_tf(t,x) + q^w(x,D_x)f(t,x)=u(t,x)\un_{\omega}(x)\,, \quad &  x \in \mathbb{R}^n, \\
f|_{t=0}=f_0 \in L^2(\rr^n),                                       &  
\end{array}\right.$$
with $q^w(x,D_x)$ being the quadratic differential operator defined by the Weyl quantization of the symbol $q$, is null-controllable from the set $\omega$ in any positive time $T>0$.
\end{theorem}

\medskip

As first examples of applications, we notice that Theorem~\ref{Main_result_control} allows us to establish the null-controllability of the harmonic heat equation
\begin{equation}\label{kk1}
\left( \partial_t - \Delta_x +|x|^2\right)f(t,x)=u(t,x)\un_{\omega}(x)\,, \quad (t,x)\in(0,T)\times\mathbb{R}^n,
\end{equation}
from any open control set satisfying condition (\ref{hyp_omega}) in any positive time. However, notice that in the one-dimensional case, this harmonic heat equation (\ref{kk1}) is shown to be not null-controllable from the half line $\omega=(a,+\infty)$ in any positive time~\cite{Duykaerts_Miller_resolvent} (Proposition~5.1). The result of Theorem~\ref{Main_result_control} also applies to the Kramers-Fokker-Planck equation (\ref{eq:KFP}):
$$\Big(\partial_t -\Delta_v+\frac{v^2}{4}+v\partial_x-ax\partial_v\Big)f(t,v,x)=u(t,v,x)\un_{\omega}(v,x)\,, \quad (t,v,x)\in(0,T)\times\mathbb{R}^2,$$
when $a \in \rr^*$.

Since the $L^2(\rr^n)$-adjoint of a quadratic operator $(q^w,D(q^w))$ is given by the quadratic operator $(\overline{q}^w,D(\overline{q}^w))$, whose Weyl symbol is the complex conjugate of the symbol $q$, we notice that the assumptions of Theorem~\ref{Main_result_control} hold for the operator $P=q^w(x,D_x)$ if and only if they hold for its $L^2(\rr^n)$-adjoint operator $P^*=\overline{q}^w(x,D_x)$. By using the Hilbert Uniqueness Method~\cite{coron_book} (Theorem~2.44), the result of null-controllability given by Theorem~\ref{Main_result_control} is therefore equivalent to the following observability estimate:

\medskip

\begin{theorem} \label{Main_result_obs}
Let $q : \rr_{x}^{n} \times \rr_{\xi}^n \rightarrow \cc$ be a complex-valued quadratic form with a non negative real part $\emph{\textrm{Re }}q \geq 0$, and a zero singular space $S=\{0\}$. If $\omega$ is an open subset of $\mathbb{R}^n$ satisfying condition (\ref{hyp_omega}), then for all $T>0$, there exists a positive constant $C_T>0$ such that 
$$\forall g \in L^2(\rr^n), \quad \| e^{-Tq^w}g\|_{L^2(\rr^n)}^2 \leq C_T\int_{0}^T\| e^{-tq^w}g\|_{L^2(\omega)}^2dt,$$
where $(e^{-tq^w})_{t \geq 0}$ denotes the contraction semigroup on $L^2(\rr^n)$ generated by the quadratic operator $q^w(x,D_x)$.
\end{theorem}

\medskip

The proof of Theorem \ref{Main_result_obs} is given in Section~\ref{sec:Proof}. It points out in particular that the control cost in the above observability estimate satisfies
$$\exists C>1, \forall T>0, \quad C_T=C\exp\Big(\frac{C}{T^{2k_0+1}}\Big),$$
where $0 \leq k_0 \leq 2n-1$ is the smallest integer verifying (\ref{h1bis2}).

We close this paragraph with a few comments on how relate the two key assumptions ensuring null-controllability in Theorems~\ref{meta_thm} and~\ref{Main_result_control}, that are respectively the Kalman rank condition (\ref{kal1}) and the condition of zero singular space $S=\{0\}$. Indeed, we notice that up to a constant the opposite of the hypoelliptic Ornstein-Uhlenbeck operator 
$$-P=-\frac{1}{2}\textrm{Tr}(Q\nabla_x^2)-\langle Bx,\nabla_x\rangle=q^w(x,D_x)+\frac{1}{2}\textrm{Tr}(B),$$ 
is a quadratic operator $q^w(x,D_x)$, whose Weyl symbol 
$$q(x,\xi)=\frac{1}{2}\langle Q\xi,\xi\rangle-i\langle Bx,\xi \rangle,$$
has a non-negative real part $\textrm{Re }q \geq 0$. A direct computation shows that the Kalman rank condition (\ref{kal1}) actually implies that the singular space of the quadratic form $q$ is equal to 
$$S=\rr_x^n\times \{0\} \subset \rr_x^{n} \times \rr_{\xi}^n.$$
The result of Theorem~\ref{meta_thm} is therefore not a byproduct of Theorem~\ref{Main_result_control}. More specifically, the condition of zero singular space $S=\{0\}$ accounts for the smoothing properties of the semigroup $(e^{-tq^w})_{t \geq 0}$ both on the direct and Fourier sides, that is, for both smoothing and decaying properties of the semigroup solution $e^{-tq^w}g$ for any positive time $t>0$. On the other hand, the condition that the singular space is equal to $S=\rr_x^n\times \{0\}$ only accounts for the smoothing properties of the semigroup solution $e^{-tq^w}g$ for any positive time $t>0$, but not for any decaying property. It explains why the proofs of Theorems~\ref{meta_thm} and~\ref{Main_result_control} rely on different smoothing properties of semigroups, namely Gevrey smoothing properties for the proof of Theorem~\ref{meta_thm}, and Gelfand-Shilov smoothing properties for the one of Theorem~\ref{Main_result_control}. It also accounts for the fact that the orthogonal projections used in the first case are frequency cutoff projections, whereas the ones used in the second case are Hermite projections.

As highlighted in~\cite{rodwahl}, the notion of singular space actually allows one to sharply understand the propagation of Gabor singularities (characterizing the lack of Schwartz regularity) of the semigroup solution $e^{-tq^w}g$ associated with any accretive quadratic operator. The lack of Schwartz regularity of a tempered distribution is characterized by its Gabor wave front set whose definition and basic properties are recalled in~\cite{rodwahl}. The Gabor wave front set (or Gabor singularities) was introduced by H\"ormander~\cite{Hormander1} and measures the directions in the phase space in which a tempered distribution does not behave like a Schwartz function. It is hence empty if and only if a distribution that is a priori tempered is in fact a Schwartz function. The Gabor wave front set thus measures global regularity in the sense of both smoothness and decay at infinity. More specifically, it is pointed out in~\cite{rodwahl} that only Gabor singularities of the initial datum $g \in L^2(\rr^n)$ contained in the singular space $S$ of the quadratic symbol $q$, can propagate for positive times along the curves given by the flow $(e^{tH_{\textrm{Im}q}})_{t \in \rr}$ of the Hamilton vector field 
$$H_{\textrm{Im}q}=\frac{\partial \textrm{Im } q}{\partial \xi} \cdot \frac{\partial}{\partial_x}-\frac{\partial \textrm{Im }q}{\partial x} \cdot \frac{\partial}{\partial_{\xi}},$$
associated with the imaginary part of the symbol. On the other hand, the Gabor singularities of the initial datum outside the singular space are all smoothed out for any positive time. More specifically, the following microlocal inclusion of Gabor wave front sets is established in~\cite{rodwahl} (Theorem~6.2),
\begin{equation}\label{etoile}
\forall g \in L^2(\rr^n), \forall  t>0, \quad WF(e^{-tq^w}g) \subset e^{tH_{\textrm{Im}q}}\big(WF(g) \cap S\big) \subset S.
\end{equation}
The microlocal inclusion (\ref{etoile}) was shown to hold as well for other types of wave front sets, as Gelfand-Shilov wave front sets~\cite{cara}, or polynomial phase space wave front sets~\cite{polywahl}, see also~\cite{karel_mehler} for a generalization of the microlocal inclusion (\ref{etoile}) of Gabor wave front sets in the non-autonomous case.

\subsubsection{Null-controllability and observability of hypoelliptic Ornstein-Uhlenbeck equations posed in weighted $L^2$-spaces with respect to invariant measures}\label{sec:LRsub1}

Let 
\begin{equation}\label{opOU}
P=\frac{1}{2}\textrm{Tr}(Q\nabla_x^2)+\langle Bx,\nabla_x\rangle, \quad x \in \rr^n,
\end{equation}
where $Q$ and $B$ are real $n \times n$-matrices, with $Q$ symmetric positive semidefinite, be a Ornstein-Uhlenbeck operator satisfying the Kalman rank condition (\ref{kal1}).

The existence of an invariant measure $\mu$ for the Markov semigroup $(e^{tP})_{t \geq 0}$ defined in (\ref{gen}), that is, a probability measure on $\rr^n$ verifying 
$$\forall t \geq 0, \forall f \in C_b(\rr^n), \quad \int_{\rr^n}(e^{tP}f)(x)d\mu(x)=\int_{\rr^n}f(x)d\mu(x),$$
where $C_b(\rr^n)$ stands for the space of continuous and bounded functions on $\rr^n$, is known to be equivalent~\cite{DaPrato} (Section~11.2.3)
to the following localization of the spectrum of~$B$,
\begin{equation}\label{pav1}
\sigma(B) \subset \mathbb{C}_-=\{z \in \mathbb{C} : \textrm{Re }z<0\}.
\end{equation}
When this condition holds, the invariant measure is unique and is given by $d\mu(x)=\rho(x)dx$, 
where the density with respect to the Lebesgue measure is
\begin{equation}\label{def:Qinfty_rho_OU}
\rho(x)=\frac{1}{(2\pi)^{\frac{n}{2}}\sqrt{\det Q_{\infty}}}e^{-\frac{1}{2}\langle Q_{\infty}^{-1}x,x\rangle}, \quad x \in \rr^n,
\end{equation}
with 
\begin{equation}\label{def:Qinfty_rho_OU1}
Q_{\infty}=\int_0^{+\infty}e^{sB}Qe^{sB^T}ds.
\end{equation}
We consider the Ornstein-Uhlenbeck operator acting on the space $L^2_{\rho}=L^2(\rr^n,\rho(x)dx)$, equipped with the domain
\begin{equation}\label{xc5}
D(P)=\{g \in L_{\rho}^2 : Pg \in L_{\rho}^2\}.
\end{equation}
On the one hand, the following result of null-controllability is an application of Theorem~\ref{Main_result_control}:

\medskip

\begin{corollary}\label{cor_Main_result_control}
Let $T>0$ and $\omega$ be an open subset of $\mathbb{R}^n$ satisfying (\ref{hyp_omega}). When the Kalman rank condition (\ref{kal1}) and the localization of the spectrum $\sigma(B) \subset \mathbb{C}_-$ hold, the Ornstein-Uhlenbeck equation posed in the $L_{\rho}^2$ space weighted by the invariant measure
\begin{equation} \label{syst_LR17}
\left\lbrace \begin{array}{ll}
\partial_t f(t,x) - \frac{1}{2}\emph{\textrm{Tr}}[Q\nabla_x^2 f(t,x)] - \langle Bx, \nabla_x f(t,x)\rangle=u(t,x)\un_{\omega}(x)\,, \quad &  x \in \mathbb{R}^n,\\
f|_{t=0}=f_0 \in L_{\rho}^2,                                       &  
\end{array}\right.
\end{equation}
is null-controllable from the set $\omega$ in time $T>0$, with a control function $u \in L^2((0,T) \times \rr^n,dt \otimes \rho(x)dx)$ supported in $[0,T]\times\omega$.
\end{corollary}

\medskip

On the other hand, the following result of observability is an application of Theorem~\ref{Main_result_obs}:

\medskip

\begin{corollary}\label{thm:OU}
Let $P$ be the Ornstein-Uhlenbeck operator defined in (\ref{opOU}),  $T>0$ and~$\omega$ be an open subset of $\mathbb{R}^n$ satisfying (\ref{hyp_omega}). When the Kalman rank condition (\ref{kal1}) and the localization of the spectrum $\sigma(B) \subset \mathbb{C}_-$ hold, then the Ornstein-Uhlenbeck operator satisfies the following observability estimate:
$$\exists C_T>0, \forall g \in L^2(\mathbb{R}^n,\rho(x) dx), \quad \| e^{TP} g\|_{L^2(\mathbb{R}^n,\rho(x) dx)}^2 \leq C_T \int_0^T \| e^{tP} g\|_{L^2(\omega,\rho(x) dx)}^2 dt,$$
where $(e^{tP})_{t \geq 0}$ denotes the semigroup on $L^2(\mathbb{R}^n,\rho(x) dx)$ generated by $-P$.
\end{corollary}

\medskip

The proofs of Corollary~\ref{cor_Main_result_control} and Corollary~\ref{thm:OU} are given in Section~\ref{sec:LRsub15}. As an application, let us mention that the results of Corollary~\ref{cor_Main_result_control} and Corollary~\ref{thm:OU} apply for instance to the system of linear stochastic differential equations given in~\cite[Section~4.2]{kps11}  that is obtained as a finite-dimensional Markovian approximation of the non-Markovian generalized Langevin equation in $\rr^n$,
\begin{equation}\label{GLE}
\ddot{x}=-\nabla_xV(x)-\int_0^t \gamma(t-s)\dot{x}(s)ds+F(t),
\end{equation}
where $V(x)=\frac{1}{2}\omega^2 x^2$ is a non-degenerate quadratic potential and $F(t)$ a mean zero stationary Gaussian process with autocorrelation function 
$$\gamma(t)=\sum_{j=1}^m\lambda_j^2e^{-\alpha_j|t|}, \qquad \alpha_j>0, \lambda_j>0,$$ 
in accordance to the fluctuation-dissipation theorem
$$\langle F(t) \otimes F(s) \rangle=\beta^{-1}\gamma(t-s)I_n, \qquad \beta>0,$$
with $I_n$ being the identity matrix. We refer the readers to the work~\cite{kps11} for further details about this model.

\subsubsection{Null-controllability and observability of hypoelliptic Fokker-Planck equations posed in weighted $L^2$-spaces with respect to invariant measures}\label{sec:LRsub23}

We consider the Fokker-Planck operator
\begin{equation}\label{jen0772}
\mathscr{P}=\frac{1}{2}\textrm{Tr}(Q\nabla_x^2)-\langle Bx,\nabla_x\rangle-\textrm{Tr}(B), \quad x \in \rr^n,
\end{equation}
where $Q$ and $B$ are real $n \times n$-matrices, with $Q$ symmetric positive semidefinite. We assume that the Kalman rank condition 
and the localization of the spectrum of~$B$,
\begin{equation}\label{kal1771}
\textrm{Rank}[Q^{\frac{1}{2}},BQ^{\frac{1}{2}},\dots, B^{n-1}Q^{\frac{1}{2}}]=n, \qquad \sigma(B) \subset \mathbb{C}_-,
\end{equation} 
hold. As before, we set
\begin{equation}\label{pav277}
\rho(x)=\frac{1}{(2\pi)^{\frac{n}{2}}\sqrt{\det Q_{\infty}}}e^{-\frac{1}{2}\langle Q_{\infty}^{-1}x,x\rangle},
\end{equation}
with 
\begin{equation}\label{pav377}
Q_{\infty}=\int_0^{+\infty}e^{sB}Qe^{sB^T}ds.
\end{equation}
We consider the operator $\mathscr{P}$ acting on the space $L^2_{1/\rho}=L^2(\rr^n,\rho(x)^{-1}dx)$, equipped with the domain
\begin{equation}\label{xc544}
D(\mathscr{P})=\{g \in L^2_{1/\rho} : \mathscr{P}g \in L^2_{1/\rho}\}.
\end{equation}
On the one hand, the following result of null-controllability is an application of Theorem~\ref{Main_result_control}:

\medskip

\begin{corollary}\label{cor_Main_result_control1}
Let $T>0$ and $\omega$ be an open subset of $\mathbb{R}^n$ satisfying (\ref{hyp_omega}). When the Kalman rank condition (\ref{kal1}) and the localization of the spectrum $\sigma(B) \subset \mathbb{C}_-$ hold, the Fokker-Planck equation posed in the $L^2_{1/\rho}$ space weighted by the invariant measure
\begin{equation} \label{syst_LR179}
\left\lbrace \begin{array}{l}
\partial_t f(t,x) - \frac{1}{2}\emph{\textrm{Tr}}[Q\nabla_x^2 f(t,x)]+\langle Bx, \nabla_x f(t,x)\rangle+\emph{\textrm{Tr}}(B)f(t,x)=u(t,x)\un_{\omega}(x),\\
f|_{t=0}=f_0 \in L^2_{1/\rho},                                       
\end{array}\right.
\end{equation}
is null-controllable from the set $\omega$ in time $T>0$, with a control function $u \in L^2((0,T) \times \rr^n,dt \otimes \rho(x)^{-1}dx)$ supported in $[0,T]\times\omega$.
\end{corollary}

\medskip

On the other hand, the following result of observability is an application of Theorem~\ref{Main_result_obs}:

\medskip

\begin{corollary}\label{thm:OU1}
Let $\mathscr{P}$ be the Fokker-Planck operator defined in (\ref{jen0772}),  $T>0$ and $\omega$ be an open subset of $\mathbb{R}^n$ satisfying (\ref{hyp_omega}). When the Kalman rank condition (\ref{kal1}) and the localization of the spectrum $\sigma(B) \subset \mathbb{C}_-$ hold, then the Fokker-Planck operator satisfies the following observability estimate:
$$\exists C_T>0, \forall g \in L^2_{1/\rho}, \quad \| e^{T\mathscr{P}} g\|_{L^2(\mathbb{R}^n,\rho^{-1}(x) dx)}^2 \leq C_T \int_0^T \| e^{t\mathscr{P}} g\|_{L^2(\omega,\rho^{-1}(x) dx)}^2 dt,$$
where $(e^{t\mathscr{P}})_{t \geq 0}$ denotes the semigroup on $L^2(\rr^n,\rho(x)^{-1}dx)$ generated by $-\mathscr{P}$.
\end{corollary}

\medskip

The proofs of Corollaries~\ref{cor_Main_result_control1} and~\ref{thm:OU1} are given in Section~\ref{sec:LRsub151}.

\subsubsection{Outline of the work} 
In Section~\ref{sec:Meta_thm}, we state a general observability estimate whose proof is given in Appendix (Section~\ref{Appendix_LM}). This proof relies on an adapted Lebeau-Robbiano method in which projection operators do not necessarily commute with semigroups.
Thanks to this general result (Theorem~\ref{Meta_thm_AdaptedLRmethod}), Theorems~\ref{meta_thm}, \ref{Main_result_control} and \ref{Main_result_obs} are then derived in an unified way.
Section~\ref{sec:LR} is devoted to the proof of null-controllability for hypoelliptic Ornstein-Uhlenbeck equations posed in the $L^2(\rr^n)$ space, whereas the proof of observability for parabolic equations associated with accretive quadratic operators with zero singular spaces is given in Section~\ref{sec:Proof}. Sections~\ref{sec:LRsub15} and~\ref{sec:LRsub151} are then devoted to the proofs of null-controllability and observability for respectively hypoelliptic Ornstein-Uhlenbeck and Fokker-Planck equations posed in $L^2$-spaces weighted by invariant measures. Section~\ref{sec:appli} provides an application of Theorems~\ref{Main_result_control} and~\ref{Main_result_obs} for the study of a model of a two oscillators chain coupled with two heat baths at each side. Section~\ref{appendix} is an appendix giving the proof of a spectral inequality for Hermite functions used in Section~\ref{sec:Proof}, a reminder about the Gelfand-Shilov regularity and the proof of the general observability estimate (Theorem~\ref{Meta_thm_AdaptedLRmethod}) written in collaboration with Luc Miller\footnote{Universit\'e Paris-Ouest, Nanterre La D\'efense, UFR SEGMI, B\^atiment G, 200 Av. de la R\'epublique, 92001 Nanterre Cedex, France (luc.miller@math.cnrs.fr)}.

\bigskip

\noindent
\textbf{Acknowledgements.} The authors are most grateful to Luc Miller and the referees for indicating missing references, 
and their very enriching remarks and stimulating comments which have helped to nicely simplify some parts of the proofs contained in this work.

\section{Adapted Lebeau-Robbiano method for observability} \label{sec:Meta_thm}

This section is devoted to the statement of the following general observability estimate,
that will allow to prove Theorems \ref{meta_thm}, \ref{Main_result_control} and \ref{Main_result_obs} in a unified way.

\medskip

\begin{theorem} \label{Meta_thm_AdaptedLRmethod}
Let $\Omega$ be an open subset of $\mathbb{R}^n$,
 $\omega$ be an open subset of $\Omega$,
 $(\pi_k)_{k \in \mathbb{N}^*}$ be a family of orthogonal projections defined on $L^2(\Omega)$,
 $(e^{tA})_{t \geq 0}$ be a contraction semigroup on $L^2(\Omega)$;
 $c_1, c_2, a, b, t_0, m >0$ be positive constants with $a<b$.
If the following spectral inequality
\begin{equation} \label{Meta_thm_IS}
\forall g \in L^2(\Omega), \forall k \geq 1, \quad \|\pi_k g \|_{L^2(\Omega)} \leq e^{c_1 k^a} \|\pi_k g \|_{L^2(\omega)},
\end{equation}
and the following dissipation estimate 
\begin{equation} \label{Meta_thm_dissip}
\forall g \in L^2(\Omega), \forall k \geq 1, \forall 0<t<t_0, \quad \| (1-\pi_k)(e^{tA} g)\|_{L^2(\Omega)} \leq \frac{1}{c_2} e^{-c_2 t^m k^b} \|g\|_{L^2(\Omega)},
\end{equation}
hold, then there exists a positive constant $C>1$ such that the following observability estimate holds
\begin{equation} \label{meta_thm_IO}
\forall T>0, \forall g \in L^2(\Omega), \quad \| e^{TA} g \|_{L^2(\Omega)}^2 \leq C\exp\Big(\frac{C}{T^{\frac{am}{b-a}}}\Big) \int_0^T \|e^{tA} g \|_{L^2(\omega)}^2 dt.
\end{equation}
\end{theorem}

\medskip

We stress the fact that the assumptions in the above statement do not require that the orthogonal projections $(\pi_k)_{k \geq 1}$ are spectral projections onto the eigenspaces of the operator $A$, which is allowed to be non-selfadjoint. We shall see in the proof that the possible lack of commutation between the contraction semigroup $(e^{tA})_{t \geq 0}$ and the orthogonal projections $(\pi_k)_{k \geq 1}$ can be compensated by the dissipation estimate (\ref{Meta_thm_dissip}).

The first version of the present work (arXiv:1603.05367) did not contain Theorem~\ref{Meta_thm_AdaptedLRmethod}, 
which is an abstract observability result of independent interest,
nor the estimate on the control cost in (\ref{meta_thm_IO}) as $T$ tends to zero.
The proof of Theorem \ref{Meta_thm_AdaptedLRmethod} is written in Appendix (Section~\ref{Appendix_LM}) in collaboration with Luc Miller.
It is inspired from the works~\cite{miller2010,miller_SMF} with a modification suggested to us by the author.
This strategy is simpler and more elegant than the one developped in the initial version of this article.
Notice that the constant $C>1$ appearing in the control cost in (\ref{meta_thm_IO}) can be expressed in terms of other rates $c_1$, $c_2$ and the exponents $a$,
$b$ and $m$ following the same optimization procedure as the one used in~\cite{miller2010}.

\section{Proof of null-controllability and observability of hypoelliptic Ornstein-Uhlenbeck equations}\label{sec:LR}

This section is devoted to the proof of Theorem~\ref{meta_thm}. 
By using the changes of unknowns $f=e^{-\frac{1}{2}\textrm{Tr}(B)t}\tilde{f}$ and $u=e^{-\frac{1}{2}\textrm{Tr}(B)t}\tilde{u}$, where $f$ is a solution to (\ref{syst_LR1}) with control $u$,
we begin by noticing that the result of Theorem~\ref{meta_thm} is equivalent to the null-controllability of the equation 
\begin{equation} \label{syst_LR1.z}
\left\lbrace \begin{array}{l}
\partial_t \tilde{f} - \frac{1}{2}\textrm{Tr}[Q\nabla_x^2 \tilde{f}] - \langle Bx, \nabla_x \tilde{f}\rangle-\frac{1}{2}\textrm{Tr}(B)\tilde{f}=\tilde{u}(t,x)\un_{\omega}(x)\,, \\
\tilde{f}|_{t=0}=f_0 \in L^2(\rr^n),                                        
\end{array}\right.
\end{equation}
from the set $\omega$ in time $T>0$, where $\omega$ is an open subset of $\mathbb{R}^n$ satisfying (\ref{hyp_omega}). We observe that the $L^2(\rr^n)$-adjoint of the operator 
$$\frac{1}{2}\textrm{Tr}(Q\nabla_x^2)+\langle Bx, \nabla_x\rangle+\frac{1}{2}\textrm{Tr}(B),$$ 
is given by
$$\Big(\frac{1}{2}\textrm{Tr}(Q\nabla_x^2)+\langle Bx, \nabla_x\rangle+\frac{1}{2}\textrm{Tr}(B)\Big)^*=\frac{1}{2}\textrm{Tr}(Q\nabla_x^2)+\langle (-B)x, \nabla_x\rangle+\frac{1}{2}\textrm{Tr}(-B).$$
By using the Hilbert Uniqueness Method~\cite{coron_book} (Theorem~2.44), the result of null-controllability of the equation (\ref{syst_LR1.z}) is equivalent to the following observability estimate
\begin{equation}\label{dfg2}
\forall T>0, \exists C_T>0, \forall g \in L^2(\rr^n), \quad \| e^{T\tilde{P}}g\|_{L^2(\rr^n)}^2 \leq C_T\int_{0}^T\| e^{t\tilde{P}}g\|_{L^2(\omega)}^2dt,
\end{equation}
with 
$$\tilde{P}=\frac{1}{2}\textrm{Tr}(Q\nabla_x^2)+\langle (-B)x, \nabla_x\rangle+\frac{1}{2}\textrm{Tr}(-B).$$
As the assumptions of Theorem~\ref{meta_thm} are fulfilled when interchanging $B$ and $-B$, it is therefore equivalent to prove the observability estimate (\ref{dfg2}) for the operator 
\begin{equation}\label{dfg1}
\tilde{P}=\frac{1}{2}\textrm{Tr}(Q\nabla_x^2)+\langle Bx, \nabla_x\rangle+\frac{1}{2}\textrm{Tr}(B),
\end{equation}
with $Q$ and $B$ real $n \times n$ matrices satisfying the assumptions of Theorem~\ref{meta_thm}. 
On the other hand, we notice that the operator $\tilde{P}=-q^w(x,D_x)$ writes as the Weyl quantization of the quadratic symbol 
$$q(x,\xi)=\frac{1}{2}\langle Q^{1/2}\xi,Q^{1/2}\xi \rangle-i\langle Bx,\xi \rangle,$$
with $\langle \cdot,\cdot \rangle$ being the Euclidean scalar product on $\rr^n$, whose real part is non-negative. As recalled above (see also e.g.~\cite{mehler}), 
the operator $-\tilde{P}$ generates a contraction semigroup on $L^2(\rr^n)$.
\\

We establish the observability inequality (\ref{dfg2}) for the operator $\widetilde{P}$ defined in (\ref{dfg1}) by applying Theorem \ref{Meta_thm_AdaptedLRmethod}.
To that end, we introduce $\pi_j:L^2(\mathbb{R}^n) \rightarrow E_j$ the orthogonal frequency cutoff projection onto the closed subspace
\begin{equation}\label{sdf6}
E_j=\big\{ f \in L^2(\mathbb{R}^n) : \text{supp}(\hat{f}) \subset \{\xi \in \rr^n : |\xi| \leq j\}\big\}, \quad j \geq 1,
\end{equation}
$|\cdot|$ being the Euclidian norm on $\mathbb{R}^n$.
The following two subsections are devoted to the proofs of a spectral inequality of type  (\ref{Meta_thm_IS}) and a dissipation estimate of type (\ref{Meta_thm_dissip}).

\subsection{Dissipation estimate}

In order to derive an explicit decay rate for the Fourier transform of the contraction semigroup $(e^{t\tilde{P}}g)_{t \geq 0}$, we need the following algebraic result:

\medskip

\begin{lemma}\label{alg}
Let $Q$ and $B$ be real $n \times n$-matrices, with $Q$ symmetric positive semidefinite. When the Kalman rank condition holds
$$\emph{\textrm{Rank}}[Q^{\frac{1}{2}},BQ^{\frac{1}{2}},\dots, B^{n-1}Q^{\frac{1}{2}}]=n,$$ 
then there exist positive constants $c>0$ and $0<t_0 \leq 1$ such that  
$$\forall 0 \leq t \leq t_0, \forall X \in \rr^{n}, \quad \int_0^t|Q^{\frac{1}{2}}e^{sB^T}X|^2ds \geq ct^{2k_0+1}|X|^2,$$
with $|\cdot|$ being the Euclidean norm on $\rr^{n}$, where $0 \leq k_0 \leq n-1$ denotes the smallest integer satisfying 
$$\emph{\textrm{Rank}}[Q^{\frac{1}{2}},BQ^{\frac{1}{2}},\dots, B^{k_0}Q^{\frac{1}{2}}]=n.$$
\end{lemma}

\medskip 

\begin{proof}
We consider the function
$$f_X(t)=\int_0^t|Q^{\frac{1}{2}}e^{sB^T}X|^2ds, \quad t \in \rr,$$
depending on the parameter $X \in \rr^{n}$.
We easily check by the Leibniz formula that 
\begin{equation}\label{re1}
\forall n \geq 0, \forall t \in \rr, \quad f_X^{(n+1)}(t)=\sum_{k=0}^n\binom{n}{k} \langle Q^{\frac{1}{2}}(B^T)^{n-k}e^{tB^T}X,Q^{\frac{1}{2}}(B^T)^{k}e^{tB^T}X \rangle,
\end{equation}
where $\langle \cdot,\cdot \rangle$ denotes the Euclidean scalar product on $\rr^{n}$. According to the Kalman rank condition, we can consider $0 \leq k_0 \leq n-1$ the smallest integer satisfying 
$$\textrm{Rank}[Q^{\frac{1}{2}},BQ^{\frac{1}{2}},\dots, B^{k_0}Q^{\frac{1}{2}}]=n.$$
We therefore have 
$$\textrm{Ran}(Q^{\frac{1}{2}})+\textrm{Ran}(BQ^{\frac{1}{2}})+...+\textrm{Ran}(B^{k_0}Q^{\frac{1}{2}})=\rr^{n}.$$
This implies that 
\begin{equation}\label{re2}
\bigcap_{j=0}^{k_0}\textrm{Ker}\big(Q^{\frac{1}{2}}(B^T)^j\big) \cap \rr^{n}=\{0\}.
\end{equation}
By induction, we easily check from (\ref{re1}) that for all $k \geq 0$, 
\begin{equation}\label{re3}
\forall 0 \leq l \leq 2k+1, \ f_X^{(l)}(0)=0 \Longleftrightarrow  X \in \bigcap_{j=0}^{k}\textrm{Ker}\big(Q^{\frac{1}{2}}(B^T)^j\big) \cap \rr^{n}.
\end{equation}
According to (\ref{re1}), (\ref{re2}) and (\ref{re3}), it follows that for all $X \in \rr^{n} \setminus \{0\}$, there exists $0 \leq \tilde{k}_X \leq k_0$ such that 
\begin{equation}\label{re0}
\forall 0 \leq j \leq 2\tilde{k}_X, \quad f_X^{(j)}(0)=0, \qquad f_X^{(2\tilde{k}_X+1)}(0)=\binom{2\tilde{k}_X}{\tilde{k}_X}|Q^{\frac{1}{2}}(B^T)^{\tilde{k}_X}X|^2>0.
\end{equation}
We aim at proving that for all $X \in \mathbb{S}^{n-1}$ in the unit sphere, there exist some positive constants $c_X>0$, $0<t_X \leq 1$ and an open neighborhood $V_X$ of $X$ in $\mathbb{S}^{n-1}$ such that
\begin{equation}\label{re4}
\forall Y \in V_{X}, \forall 0 \leq t \leq t_X, \quad \int_0^t|Q^{\frac{1}{2}}e^{sB^T}Y|^2ds \geq c_X  t^{2\tilde{k}_X+1}.
\end{equation}
By analogy with~\cite[Proposition~3.2]{joh}, we proceed by contradiction. If the assertion (\ref{re4}) does not hold, there exist a sequence of positive real numbers $(t_l)_{l \geq 0}$ and a sequence $(Y_l)_{l \geq 0}$ of elements in $\mathbb{S}^{n-1}$ so that 
\begin{equation}\label{re5}
\lim_{l \to +\infty}t_l=0, \quad \lim_{l \to +\infty}Y_l=X, \quad \textrm{ and } \quad  \lim_{l \to +\infty}\frac{1}{t_l^{2\tilde{k}_X+1}}\int_0^{t_l}|Q^{\frac{1}{2}}e^{sB^T}Y_l|^2ds=0.
\end{equation}
We deduce from (\ref{re5}) that 
\begin{equation}\label{re6}
\lim_{l \to +\infty}\frac{1}{t_l^{2\tilde{k}_X+1}}\sup_{0 \leq t \leq t_l}\int_0^{t}|Q^{\frac{1}{2}}e^{sB^T}Y_l|^2ds=0.
\end{equation}
Setting
\begin{equation}\label{re7}
u_l(x)=\frac{1}{t_l^{2\tilde{k}_X+1}}\int_0^{xt_l}|Q^{\frac{1}{2}}e^{sB^T}Y_l|^2ds \geq 0, \quad 0 \leq x \leq 1,
\end{equation}
we can reformulate (\ref{re6}) as 
\begin{equation}\label{re8}
\lim_{l \to +\infty}\sup_{0 \leq x \leq 1}u_l(x)=0.
\end{equation}
By writing that
$$f_{Y_l}(t)=\int_0^{t}|Q^{\frac{1}{2}}e^{sB^T}Y_l|^2ds=\sum_{j=0}^{2\tilde{k}_X+1}a_l^{(j)}t^{j}+\mathcal{O}(t^{2\tilde{k}_X+2}),$$ 
when $t \to 0$, with $a_l^{(j)}=f_{Y_l}^{(j)}(0)(j!)^{-1}$, where the term $\mathcal{O}(t^{2\tilde{k}_X+2})$ appearing in the right-hand-side of the above formula can be assumed  to be independent on the integer $l$ thanks to Taylor formula with integral remainder and the fact that $(Y_l)_{l \geq 0}$ are elements of the unit sphere $\mathbb{S}^{n-1}$, we notice that 
\begin{equation}\label{re9}
u_l(x)=\sum_{j=0}^{2\tilde{k}_X+1}\frac{a_l^{(j)}}{t_l^{2\tilde{k}_X+1-j}}x^{j}+\mathcal{O}(t_l x^{2\tilde{k}_X+2}).
\end{equation}
It follows from (\ref{re5}), (\ref{re8}) and (\ref{re9}) that 
\begin{equation}\label{re10}
\lim_{l \to +\infty}\sup_{0 \leq x \leq 1}|p_l(x)|=0,
\end{equation}
with 
\begin{equation}\label{re11}
p_l(x)=\sum_{j=0}^{2\tilde{k}_X+1}\frac{a_l^{(j)}}{t_l^{2\tilde{k}_X+1-j}}x^{j}.
\end{equation}
By using the equivalence of norms in finite-dimensional vector space, we deduce from (\ref{re10}) that 
\begin{equation}\label{re12}
\forall 0 \leq j \leq 2\tilde{k}_X+1, \quad \lim_{l \to +\infty}\frac{a_l^{(j)}}{t_l^{2\tilde{k}_X+1-j}}=0.
\end{equation}
We obtain in particular that 
\begin{equation}\label{re13}
\lim_{l \to +\infty}a_l^{(2\tilde{k}_X+1)}=0.
\end{equation}
According to (\ref{re0}), this is in contradiction with the fact that 
\begin{equation}\label{re14}
\lim_{l \to +\infty}a_l^{(2\tilde{k}_X+1)}=\lim_{l \to +\infty}\frac{f_{Y_l}^{(2\tilde{k}_X+1)}(0)}{(2\tilde{k}_X+1)!}=\frac{f_X^{(2\tilde{k}_X+1)}(0)}{(2\tilde{k}_X+1)!}>0.
\end{equation}
By covering the compact set $\mathbb{S}^{n-1}$ by finitely many open neighborhoods of the form $(V_{X_j})_{1 \leq j \leq N}$, and letting $c=\inf_{1 \leq j \leq N}c_{X_j}>0$, $0<t_0=\inf_{1 \leq j \leq N}t_{X_j} \leq 1$, we conclude that 
$$\forall X \in \rr^n, \forall 0 \leq t \leq t_0, \quad \int_0^t|Q^{\frac{1}{2}}e^{sB^T}X|^2ds \geq c  t^{2k_0+1}|X|^2,$$
since $0 \leq \tilde{k}_X \leq k_0$. This ends the proof of Lemma~\ref{alg}.
\end{proof}

\medskip

We prove the following dissipation estimate:

\medskip

\begin{proposition} \label{thm:GS.z}
When the Kalman rank condition (\ref{kal1}) holds, then we have 
\begin{multline}\label{eq6.z}
\forall T>0, \exists C_T>1, \forall 0 \leq t \leq T, \forall k \geq 0, \forall g_0 \in L^2(\rr^n), \\ \| (1-\pi_k)(e^{t\tilde{P}}g_0)\|_{L^2(\rr^n)} \leq e^{- \delta(t)k^2} \|g_0\|_{L^2(\rr^n)},
\end{multline} 
with $0 \leq k_0 \leq n-1$ being the smallest integer satisfying 
$$\emph{\textrm{Rank}}[Q^{\frac{1}{2}},BQ^{\frac{1}{2}},\dots, B^{k_0}Q^{\frac{1}{2}}]=n$$ and
\begin{equation}\label{eq5.z}
\delta(t)=\frac{1}{C_T}\inf(t,t_0)^{2k_0+1} \geq 0, \quad t \geq 0,
\end{equation} 
with $0<t_0 \leq 1$ being defined in Lemma~\ref{alg}.
\end{proposition}

\medskip

\begin{proof}
Let $g_0 \in L^2(\mathbb{R}^n)$ and $g(t)=e^{t\tilde{P}}g_0$ be the solution of
$$\left\lbrace \begin{array}{l}
\partial_t g(t,x) - \frac{1}{2}\textrm{Tr}[Q\nabla_x^2 g(t,x)] - \langle Bx , \nabla_x g(t,x) \rangle -\frac{1}{2}\textrm{Tr}(B)g(t,x)= 0\,,\\
g(0,x)=g_0(x).
\end{array}\right.$$
Then, the function $h$ uniquely defined by $g(t,x)=h(t,e^{tB}x)e^{\frac{1}{2}\textrm{Tr}(B)t}$ solves
$$\left\lbrace \begin{array}{ll}
\partial_t h(t,y) - \frac{1}{2}\textrm{Tr}[e^{tB}Q e^{tB^T}\nabla_y^2  h(t,y)]  = 0\,, \quad & (t,y)\in(0,+\infty)\times\mathbb{R}^n\,,\\
h(0,y)=g_0(y)\,, & y \in \mathbb{R}^n.
\end{array}\right.$$
Thus, we obtain that for all $(t,\xi)\in[0,+\infty)\times\mathbb{R}^n$,
$$\widehat{h}(t,\xi) =  \widehat{g}_0(\xi) e^{-\frac{1}{2}\int_0^t | Q^{1/2}e^{sB^T} \xi|^2 ds},$$
implying that the Fourier transform of the semigroup $g(t)=e^{t\tilde{P}}g_0$ is given by
\begin{multline}\label{sdf1}
\widehat{g}(t,\xi) = |\text{det}(e^{-tB})|\widehat{h}( t ,e^{-tB^T} \xi)e^{\frac{1}{2}\textrm{Tr}(B)t}\\
=e^{-\frac{1}{2}t \textrm{Tr}(B)}\widehat{g_0}(e^{-tB^T} \xi)e^{-\frac{1}{2}\int_0^t | Q^{1/2}e^{(s-t)B^T} \xi|^2 ds}.
\end{multline}
We deduce from (\ref{sdf1}) and Lemma~\ref{alg} that for all $0 \leq t \leq T$, $k \geq 0$, $g_0 \in L^2(\rr^n)$,
\begin{multline}\label{sdf3}
\| (1-\pi_k)(e^{t\tilde{P}}g_0)\|_{L^2(\rr^n)}^2=\frac{e^{-t \textrm{Tr}(B)}}{(2\pi)^n}\int_{|\xi| \geq k}|\widehat{g_0}(e^{-tB^T} \xi)|^2e^{-\int_0^t | Q^{1/2}e^{(s-t)B^T} \xi|^2 ds}d\xi\\
=\frac{1}{(2\pi)^n}\int_{|e^{tB^T}\xi| \geq k}|\widehat{g_0}(\xi)|^2e^{-\int_0^t | Q^{1/2}e^{sB^T} \xi|^2 ds}d\xi
\leq \frac{1}{(2\pi)^n}\int_{|e^{tB^T}\xi| \geq k}|\widehat{g_0}(\xi)|^2e^{-\tilde{\delta}(t)|\xi|^2}d\xi,
\end{multline}
with $\tilde{\delta}(t)=c\inf(t,t_0)^{2k_0+1}$. It follows from (\ref{sdf3}) that for all $0 \leq t \leq T$, $k \geq 0$, $g_0 \in L^2(\rr^n)$,
\begin{multline}\label{sdf4}
\| (1-\pi_k)(e^{t\tilde{P}}g_0)\|_{L^2(\rr^n)}^2
\leq \frac{1}{(2\pi)^n}\int_{|\xi| \geq k e^{-t\|B\|}}|\widehat{g_0}(\xi)|^2e^{-\tilde{\delta}(t)|\xi|^2}d\xi\\
\leq e^{-\tilde{\delta}(t)k^2e^{-2t\|B\|}}\|g_0\|_{L^2(\rr^n)}^2\leq e^{-\tilde{\delta}(t)k^2e^{-2T\|B\|}}\|g_0\|_{L^2(\rr^n)}^2.
\end{multline}
We deduce from (\ref{sdf4}) that there exists $C_T>1$ such that for all $0 \leq t \leq T$, $k \geq 0$, $g_0 \in L^2(\rr^n)$,
\begin{equation}\label{sdf5}
\| (1-\pi_k)(e^{t\tilde{P}}g_0)\|_{L^2(\rr^n)}^2 \leq e^{-2\delta(t)k^2}\|g_0\|_{L^2(\rr^n)}^2,
\end{equation}
with 
$$\delta(t)=\frac{1}{C_T}\inf(t,t_0)^{2k_0+1} \geq 0, \quad t \geq 0.$$
It proves the estimate (\ref{eq6.z}) and ends the proof of Proposition~\ref{thm:GS.z}.
\end{proof}

\subsection{Spectral inequality for Fourier modes}

The following spectral inequality is proved by Le Rousseau and Moyano in \cite[Theorem 3.1]{LeRousseau_Moyano}:

\medskip

\begin{theorem} \label{thm:IS_LR_M}
If $\omega$ is an open subset of $\mathbb{R}^n$ satisfying condition (\ref{hyp_omega}), then there exists a positive constant $c_1>1$ such that 
$$\|\pi_k g\|_{L^2(\mathbb{R}^n)} \leq e^{c_1 k} \|\pi_k g\|_{L^2(\omega)},$$
for all $k \in \mathbb{N}^*$ and $g \in L^2(\mathbb{R}^n)$.
\end{theorem}


\subsection{Proof of Theorem~\ref{meta_thm}}
We deduce from Theorem~\ref{Meta_thm_AdaptedLRmethod} with the following choices of parameters:
\begin{itemize}
\item[$(i)$] $\Omega=\mathbb{R}^n$,
\item[$(ii)$]  $A=\tilde{P}$,
\item[$(iii)$]  $a=1$, $b=2$, 
\item[$(iv)$]  $m=2k_0+1$, where $k_0$ is defined in Proposition~\ref{thm:GS.z},
\item[$(v)$]  $t_0>0$ as in Proposition \ref{thm:GS.z},
\item[$(vi)$]  $0<c_2=1/C_{t_0}<1$, where $C_{t_0}>1$ is the constant defined in Proposition~\ref{thm:GS.z},
\item[$(vii)$]  $c_1>0$ as in Theorem \ref{thm:IS_LR_M},
\end{itemize}
that 
$$\exists C>1, \forall T>0, \forall g \in L^2(\rr^n), \quad \| e^{T\tilde{P}} g \|_{L^2(\rr^n)}^2 \leq C\exp\Big(\frac{C}{T^{2k_0+1}}\Big) \int_0^T \|e^{t\tilde{P}} g \|_{L^2(\omega)}^2 dt.$$
It proves the observability estimate (\ref{dfg2}) and ends the proof of Theorem~\ref{meta_thm}.

\section{Proof of null-controllability and observability of parabolic equations associated with accretive quadratic operators with zero singular spaces}\label{sec:Proof}

This section is devoted to the proof of Theorem~\ref{Main_result_obs}. 
As in the previous section, we use the general observability estimate established in Theorem~\ref{Meta_thm_AdaptedLRmethod}. 
Indeed, the classical Lebeau-Robbiano method cannot be directly applied in its usual form as the (generalized) eigenfunctions 
of accretive quadratic operators with zero singular spaces do not constitute in general a $L^2$-Hilbert basis. 
Contrary to the usual Lebeau-Robbiano strategy, the solutions are therefore not expanded on the (generalized) eigenfunctions 
of the operator defining the parabolic equation, but in the $L^2$-Hermite basis that does not diagonalize the operator. 
With this choice, the difficulty is that the semigroup is not diagonal anymore in the $L^2$-Hermite basis, 
and even if a finite number of modes could be steered to zero at some time, any passive control phase in the Lebeau-Robbiano method makes them all revive again. 
To overcome this lack of commutation between semigroups and Hermite projection operators, 
we take a key advantage of the Gelfand-Shilov regularizing properties of semigroups generated by accretive quadratic operators with zero singular spaces, 
and the fact that Gelfand-Shilov regularity is characterized by a certain exponential decay of the Hermite coefficients.

\subsection{Gelfand-Shilov regularizing properties}

In the following, we denote 
\begin{equation}\label{estt3}
\mathbb{P}_k g =\sum_{\substack{\alpha \in \N^n\\ |\alpha|=k}}(g,\psi_{\alpha})_{L^2(\rr^n)}\psi_{\alpha}, \quad k \geq 0, 
\end{equation} 
the orthogonal projection onto the $k^{\textrm{th}}$ energy level associated with the harmonic oscillator 
$$\mathcal{H}=-\Delta_x+|x|^2=\sum_{k=0}^{+\infty}(2k+n)\mathbb{P}_k,$$ 
where $(\psi_{\alpha})_{\alpha \in \N^n}$ stands for the $L^2$-Hermite basis. We also consider the orthogonal projection
\begin{equation}\label{estt4}
\pi_k=\sum_{j=0}^k\mathbb{P}_j, \quad k \geq 0,
\end{equation}
onto energy levels less than or equal to $k$. The exponential decay results given by the following proposition are key byproducts of the Gelfand-Shilov regularizing properties of semigroups generated by accretive quadratic operators with zero singular spaces:

\medskip

\begin{proposition} \label{thm:GS}
Let $q : \rr_{x,\xi}^{2n} \rightarrow \cc$ be a quadratic form with a non-negative real part $\emph{\textrm{Re }}q \geq 0$ and a zero singular space $S=\{0\}$. 
There exist some positive constants $C_0>1$ and $t_0>0$ such that for all $t \geq 0$, $k \geq 0$, $g \in L^2(\rr^n)$, 
\begin{equation}\label{eq6}
\| (1-\pi_k)(e^{-tq^w}g)\|_{L^2(\rr^n)} \leq C_0e^{- \delta(t)k} \|g\|_{L^2(\rr^n)},
\end{equation} 
with $0 \leq k_0 \leq 2n-1$ being the smallest integer satisfying (\ref{h1bis2}) and
\begin{equation}\label{eq5}
\delta(t)=\frac{\inf(t,t_0)^{2k_0+1}}{C_0} \geq 0, \quad t \geq 0.
\end{equation} 
\end{proposition}

\medskip

\begin{proof}
Let $q : \rr_{x,\xi}^{2n} \rightarrow \cc$ be a quadratic form with a non-negative real part $\textrm{Re }q \geq 0$ and a zero singular space $S=\{0\}$. We recall from~\cite[p. 426]{mehler} that the quadratic operator $q^w(x,D_x)$ obtained by the Weyl quantization of the symbol $q$ is accretive and generates a contraction semigroup on $L^2(\rr^n)$. 
We denote $0 \leq k_0 \leq 2n-1$ the smallest integer satisfying (\ref{h1bis2}). In the work~\cite[Theorem~1.2]{HPSVII}, Hitrik, Viola and the second author have shown that the contraction semigroup $(e^{-tq^w})_{t \geq 0}$ is smoothing for any positive time in the Gelfand-Shilov space $S_{1/2}^{1/2}(\rr^n)$,
$$\forall g \in L^2(\rr^n), \forall t>0, \quad e^{-tq^w}g \in S_{1/2}^{1/2}(\rr^n).$$
We refer the reader to the appendix (Section~\ref{GSreg}) for the definition and some characterizations of the Gelfand-Shilov regularity. 
More specifically, we deduce from~\cite[Proposition~4.1]{HPSVII} that there exist some positive constants $C_0>1$ and $t_0>0$ such that
\begin{equation}\label{eq3}
\forall 0 \leq t \leq t_0, \quad \Big\| e^{\frac{t^{2k_0+1}}{C_0}(-\Delta_x+|x|^2)}e^{-tq^w} \Big\|_{\mathcal{L}(L^2(\rr^n))} \leq C_0,
\end{equation} 
with $\mathcal{L}(L^2(\rr^n))$ the space of bounded operators on $L^2(\rr^n)$, that is,
\begin{equation}\label{eq3..1}
\forall 0 \leq t \leq t_0, \forall g \in L^2(\rr^n), \ \sum_{\alpha \in \nn^n}|(e^{-tq^w}g,\psi_{\alpha})_{L^2(\rr^n)}|^2e^{\frac{t^{2k_0+1}}{C_0}(4|\alpha|+2n)} \leq C_0^2 \|g\|_{L^2(\rr^n)}^2.
\end{equation} 
We obtain from (\ref{eq3}) and the contraction semigroup property satisfied by $(e^{-t q^w})_{t \geq 0}$ that 
\begin{equation}\label{eq4}
\forall t \geq 0, \quad \big\| e^{\delta(t)(-\Delta_x+|x|^2)}e^{-tq^w} \big\|_{\mathcal{L}(L^2(\rr^n))} \leq C_0.
\end{equation} 
It follows from (\ref{eq4}) that for all $t \geq 0$, $k \geq 0$, $g \in L^2(\rr^n)$, 
\begin{align*}
& \ \| (1-\pi_k)(e^{-tq^w}g)\|_{L^2(\rr^n)} 
= \|(1-\pi_k)(e^{-\delta(t)(-\Delta_x+|x|^2)}e^{\delta(t)(-\Delta_x+|x|^2)}e^{-tq^w}g)\|_{L^2(\rr^n)} \\ 
= & \ \| e^{-\delta(t)(-\Delta_x+|x|^2)}(1-\pi_k)(e^{\delta(t)(-\Delta_x+|x|^2)}e^{-tq^w}g)\|_{L^2(\rr^n)}\\ 
\leq & \ e^{-\delta(t)(2(k+1)+n)}\| (1-\pi_k)(e^{\delta(t)(-\Delta_x+|x|^2)}e^{-tq^w}g)\|_{L^2(\rr^n)}\\ 
\leq & \ e^{- \delta(t) k }\| e^{\delta(t)(-\Delta_x+|x|^2)}e^{-tq^w}g\|_{L^2(\rr^n)} \leq  C_0e^{- \delta(t) k }\|g\|_{L^2(\rr^n)}.
\end{align*} 
It ends the proof of Proposition~\ref{thm:GS}.
\end{proof}

\subsection{Spectral inequality for Hermite functions}

The following spectral inequality for Hermite functions is proved in Appendix (Section~\ref{refherm}):

\medskip

\begin{proposition}\label{thm:Spect_Ineq}
Let $\omega$ be an open subset of $\mathbb{R}^n$ satisfying (\ref{hyp_omega}) and $(\psi_{\alpha})_{\alpha \in \N^n}$ the Hermite basis of $L^2(\rr^n)$ diagonalizing the harmonic oscillator $\mathcal{H}=-\Delta_x+|x|^2$. 
There exists a positive constant $C_1>1$ such that for all $k \geq 0$ and $(b_\alpha)_{\alpha \in \mathbb{N}^n} \in \mathbb{C}^{\mathbb{N}^n}$,
$$\Big(\sum_{|\alpha| \leq k} |b_\alpha|^2\Big)^{1/2}=\Big(\int\limits_{\mathbb{R}^n}\Big| \sum_{|\alpha| \leq k} b_\alpha \psi_\alpha(x)\Big|^2 dx\Big)^{1/2}
\leq  C_1 e^{C_1 \sqrt{k}} \Big(\int\limits_{\omega}\Big| \sum_{|\alpha| \leq k} b_\alpha \psi_\alpha(x)\Big|^2 dx\Big)^{1/2}\,.$$
In particular, the following estimate holds
\begin{equation}\label{eq7}
\forall k \geq 0, \forall g \in L^2(\mathbb{R}^n), \quad \|\pi_k g\|_{L^2(\rr^n)} \leq C_1e^{C_1\sqrt{k}} \|\pi_k g\|_{L^2(\omega)}\,.
\end{equation}
\end{proposition}

\medskip

Notice that when $\omega$ is a bounded set (hence does not satisfy (\ref{hyp_omega})), then the weaker
spectral inequality obtained from (\ref{eq7}) while replacing $\sqrt{k}$ by $k$ fails even
for single Hermite functions instead of sums, by at least an extra logarithmic factor
in the exponentials when $n\geq 2$. This fact is proved in~\cite{Mil08} (Section~4.2). Whether the spectral inequality (\ref{eq7}) holds with $\sqrt{k}$ replaced by $k\ln k$ in the exponential term
when $\omega$ is a bounded set is still open.
In the one-dimensional case, when $\omega$ is an open half line and $0<s<1$, the weaker spectral
inequality obtained from (\ref{eq7}) with $\sqrt{k}$ replaced by $k^s$ also fails. Indeed,
the Lebeau-Robbiano strategy would otherwise allow to establish a null-controllability result that is disproved in~\cite{Mil08}
(Section~4.3.1), see also~\cite{Duykaerts_Miller_resolvent}. 
Similarly, Remark 1.9 in~\cite{Mil08} points out that if $n\geq 3$, $\omega$ is a non-empty
open cone $\Gamma =\{ x \in\R^n : |x|>r_{0}, \frac{x}{|x|} \in \Omega_{0} \}$
$r_{0}\geq 0$ and $\Omega_{0}$ is an open subset of the unit sphere, and if there exists
a vector space of dimension~$2$ in $\R^n$ not intersecting the closure of
$\Omega_0$, then the weaker spectral inequality obtained from (\ref{eq7}) while replacing
$\sqrt{k}$ by $k^s$ fails for all $1/2< s<1$.


\subsection{Proof of Theorem~\ref{Main_result_obs}}
We deduce from Theorem~\ref{Meta_thm_AdaptedLRmethod} with the following choices of parameters:
\begin{itemize}
\item[$(i)$] $\Omega=\mathbb{R}^n$,
\item[$(ii)$] $A=-q^w(x,D_x)$,
\item[$(iii)$] $a=\frac{1}{2}$, $b=1$, 
\item[$(iv)$] $t_0>0$ as in Proposition \ref{thm:GS},
\item[$(v)$] $m=2k_0+1$, where $k_0$ is defined in Proposition~\ref{thm:GS},
\item[$(vi)$] any constant $c_1>0$ satisfying for all $k \geq 1$,  $C_1 e^{C_1 \sqrt{k}} \leq e^{c_1 \sqrt{k}}$, where the constant $C_1>1$ is defined in Proposition~\ref{thm:Spect_Ineq},
\item[$(vii)$] $c_2= \frac{1}{C_0}>0$, where $C_0>1$ is defined in Proposition \ref{thm:GS},
\end{itemize}
that 
$$\exists C>1, \forall T>0, \forall g \in L^2(\rr^n), \quad \|e^{-Tq^w} g \|_{L^2(\rr^n)}^2 \leq C\exp\Big(\frac{C}{T^{2k_0+1}}\Big) \int_0^T \|e^{-tq^w} g \|_{L^2(\omega)}^2 dt.$$
It ends the proof of Theorem~\ref{Main_result_obs}.

\section{Proofs of null-controllability and observability of hypoelliptic Ornstein-Uhlenbeck equations posed in weighted $L^2$-spaces}\label{sec:LRsub15}

Let $P$ be a hypoelliptic Ornstein-Uhlenbeck operator (\ref{opOU}) such that the Kalman rank condition (\ref{kal1}) and the localization of the spectrum $\sigma(B) \subset \mathbb{C}_-$ hold.
We consider the operator $P$ acting on the space $L^2(\rr^n,\rho(x)dx)$, with $\rho$ being the density function defined in (\ref{def:Qinfty_rho_OU}). The Kalman rank condition 
$$\textrm{Rank}[Q^{\frac{1}{2}},BQ^{\frac{1}{2}},\dots, B^{n-1}Q^{\frac{1}{2}}]=n,$$
allows one to consider $0 \leq k_0 \leq n-1$ the smallest integer satisfying
$$\textrm{Rank}[Q^{\frac{1}{2}},BQ^{\frac{1}{2}},\dots, B^{k_0}Q^{\frac{1}{2}}]=n.$$
We associate to the operator $P$ acting on $L^2_{\rho}=L^2(\rr^n,\rho(x)dx)$, the quadratic operator $\mathscr{L}$ acting on $L^2=L^2(\rr^n,dx)$,
\begin{equation}\label{pav4.5}
\mathscr{L}h=-\sqrt{\rho}P\big((\sqrt{\rho})^{-1}h\big)-\frac{1}{2}\textrm{Tr}(B)h.
\end{equation}
Recalling the notation (\ref{not}), a direct computation led in the work~\cite{uhlenb} (see (3.7) in Section~3.1) shows that
\begin{equation}\label{pav5}
\mathscr{L}=\frac{1}{2}|Q^{\frac{1}{2}}D_x|^2+\frac{1}{8}|Q^{\frac{1}{2}}Q_{\infty}^{-1}x|^2-i\Big\langle\Big(\frac{1}{2}QQ_{\infty}^{-1}+B\Big)x,D_x\Big\rangle, 
\end{equation} 
with $D_x=i^{-1}\nabla_x$, where $Q_{\infty}$ is the symmetric positive definite matrix (\ref{def:Qinfty_rho_OU1}).
The operator $\mathscr{L}=q^w(x,D_x)$ is a quadratic operator whose Weyl symbol
\begin{equation}\label{pav6}
q(x,\xi)=\frac{1}{2}|Q^{\frac{1}{2}}\xi|^2+\frac{1}{8}|Q^{\frac{1}{2}}Q_{\infty}^{-1}x|^2-i\Big\langle\Big(\frac{1}{2}QQ_{\infty}^{-1}+B\Big)x,\xi\Big\rangle, \quad (x,\xi) \in \rr^{2n},
\end{equation}
has a non-negative real part $\textrm{Re }q \geq 0$.
On the other hand, we prove in~\cite{uhlenb}, see formulas (3.22), (3.23) and (3.24), that the singular space of the quadratic operator $\mathscr{L}$ is zero 
$S=\{0\}$. 
More precisely, we show in~\cite{uhlenb} that the smallest integer $0 \leq k_0 \leq 2n-1$ satisfying 
\begin{equation}\label{pav30}
\Big(\bigcap_{j=0}^{k_0}\textrm{Ker}
\big[\textrm{Re }F(\textrm{Im }F)^j \big]\Big)\cap \rr^{2n}=\{0\},
\end{equation}
with $F$ being the Hamilton map of $q$,
corresponds exactly to the smallest integer $0 \leq k_0 \leq n-1$ satisfying
\begin{equation}\label{pav31}
\textrm{Rank}[Q^{\frac{1}{2}},BQ^{\frac{1}{2}},\dots, B^{k_0}Q^{\frac{1}{2}}]=n.
\end{equation}

Let $\omega$ be an open subset of $\mathbb{R}^n$ satisfying condition (\ref{hyp_omega}).
We can therefore deduce from Theorem~\ref{Main_result_control} applied to the quadratic operator $\mathscr{L}$ that the parabolic equation 
\begin{equation}\label{evo1}
\left\lbrace\begin{array}{c}
\partial_th(t,x)+\mathscr{L}h(t,x)=u(t,x)\un_{\omega}(x),\\
h|_{t=0}=h_0 \in L^2(\rr^n,dx),
\end{array}\right.
\end{equation}
is null-controllable from the set $\omega$ in any positive time $T>0$. Let $f_0 \in L_{\rho}^2$.
By using that the mappings 
\begin{equation}\label{evo2}
\begin{array}{cc}
\mathcal{T} :  L^2_{\rho} & \rightarrow  L^2\\
\ \ \ v &  \mapsto  \sqrt{\rho}v
\end{array}, \qquad \begin{array}{cc}
\mathcal{T}^{-1} :  L^2 & \rightarrow L^2_{\rho}\\
\quad \ \ \ v & \mapsto \sqrt{\rho}^{-1}v
\end{array},
\end{equation}
are isometric, we consider a control function $u \in L^2((0,T) \times \rr^n,dt \otimes dx)$ supported in $[0,T]\times\omega$ such that the mild solution $h$ to the equation (\ref{evo1}) with initial datum $h_0=\mathcal{T}f_0$ satisfies $h(T,\cdot)=0$. We deduce from (\ref{opOU}), (\ref{pav4.5}) and (\ref{evo2}) that the mild solution
$$f=e^{-\frac{t}{2}\textrm{Tr}(B)}\mathcal{T}^{-1}h \in L_{\rho}^2,$$ 
to the equation
$$\left\lbrace \begin{array}{ll}
\partial_t f(t,x) - \frac{1}{2}\textrm{Tr}[Q\nabla_x^2 f(t,x)] - \langle Bx, \nabla_x f(t,x)\rangle=\tilde{u}(t,x)\un_{\omega}(x)\,, \quad &  x \in \mathbb{R}^n,\\
f|_{t=0}=f_0 \in L_{\rho}^2,                                       &  
\end{array}\right.$$
with the control function supported in $[0,T]\times\omega$,
$$\tilde{u}(t,x)=e^{-\frac{t}{2}\textrm{Tr}(B)}\mathcal{T}^{-1}u(t,x) \in L^2((0,T) \times \rr^n,dt \otimes \rho(x)dx),$$
satisfies $f(T,\cdot)=0$. This proves that the hypoelliptic Ornstein-Uhlenbeck equation (\ref{syst_LR17}) is null-controllable from the set $\omega$ in any positive time $T>0$. This ends the proof of Corollary~\ref{cor_Main_result_control}.

On the other hand, we deduce from Theorem~\ref{Main_result_obs} applied to the quadratic operator $\mathscr{L}$ that for all $T>0$, there exists a positive constant $C_T>0$ such that 
\begin{equation}\label{bn1.0}
\forall h_0 \in L^2(\rr^n), \quad \| e^{-T\mathscr{L}}h_0\|_{L^2(\rr^n,dx)}^2 \leq C_T\int_{0}^T \| e^{-t \mathscr{L}}h_0\|_{L^2(\omega,dx)}^2dt.
\end{equation}
According to (\ref{pav4.5}), the semigroup $(e^{tP})_{t \geq 0}$ on $L^2_{\rho}$ is given by 
\begin{equation}\label{bn1}
\forall f_0 \in L^2_{\rho}, \forall t \geq 0, \quad  e^{tP}f_0=e^{-\frac{t}{2}\textrm{Tr}(B)}\mathcal{T}^{-1}e^{-t\mathscr{L}}\mathcal{T}f_0,
\end{equation}
where $(e^{-t\mathscr{L}})_{t \geq 0}$ denotes the $L^2(\rr^n,dx)$ contraction semigroup generated by $\mathscr{L}$. Notice from the localization of the spectrum $\sigma(B) \subset \cc_-$ of  $B \in M_n(\rr)$ that $\textrm{Tr}(B) < 0$.
By observing that 
$$\|e^{tP}f_0\|_{L^2(\rr^n,\rho(x) dx)}=e^{-\frac{t}{2}\textrm{Tr}(B)}\|e^{-t\mathscr{L}}\mathcal{T}f_0\|_{L^2(\rr^n,dx)}$$
and 
$$\|e^{tP}f_0\|_{L^2(\omega,\rho(x) dx)}=e^{-\frac{t}{2}\textrm{Tr}(B)}\|e^{-t\mathscr{L}}\mathcal{T}f_0\|_{L^2(\omega,dx)},$$
we deduce from (\ref{bn1.0}) and (\ref{bn1}) that the hypoelliptic Ornstein-Uhlenbeck operator $P$ satisfies the observability estimate
$$\forall g \in L^2(\mathbb{R}^n,\rho(x) dx), \quad \| e^{TP} g\|_{L^2(\mathbb{R}^n,\rho(x) dx)}^2 \leq \tilde{C}_T \int_0^T \| e^{tP} g\|_{L^2(\omega,\rho(x) dx)}^2 dt,$$
with $\tilde{C}_T=e^{-T\textrm{Tr}(B)}C_T>0$.
This ends the proof of Corollary~\ref{thm:OU}.

\section{Proofs of null-controllability and observability of hypoelliptic Fokker-Planck equations posed in weighted $L^2$-spaces}\label{sec:LRsub151}

Let $\mathscr{P}$ be a Fokker-Planck operator (\ref{jen0772}) satisfying conditions (\ref{kal1771}). We consider the operator $\mathscr{P}$ acting on the space $L^2(\rr^n,\rho(x)^{-1}dx)$. The Kalman rank condition 
$$\textrm{Rank}[Q^{\frac{1}{2}},BQ^{\frac{1}{2}},\dots, B^{n-1}Q^{\frac{1}{2}}]=n,$$
allows one to consider $0 \leq k_0 \leq n-1$ the smallest integer satisfying
$$\textrm{Rank}[Q^{\frac{1}{2}},BQ^{\frac{1}{2}},\dots, B^{k_0}Q^{\frac{1}{2}}]=n.$$
We associate to the operator $\mathscr{P}$ acting on $L^2_{1/\rho}=L^2(\rr^n,\rho(x)^{-1}dx)$, the quadratic operator $\mathfrak{L}$ acting on $L^2=L^2(\rr^n,dx)$,
\begin{equation}\label{pav4.59}
\mathfrak{L}h=-\sqrt{\rho}^{-1}\mathscr{P}\big(\sqrt{\rho}h\big)-\frac{1}{2}\textrm{Tr}(B)h.
\end{equation}
Recalling the notation (\ref{not}), a direct computation led in the work~\cite{uhlenb} (see (2.54) in Section~2.6) shows that
\begin{equation}\label{pav59}
\mathfrak{L}=\frac{1}{2}|Q^{\frac{1}{2}}D_x|^2+\frac{1}{8}|Q^{\frac{1}{2}}Q_{\infty}^{-1}x|^2+i\Big\langle \Big(\frac{1}{2} QQ_{\infty}^{-1}+B\Big)x,D_x\Big\rangle, 
\end{equation} 
with $D_x=i^{-1}\nabla_x$, where $Q_{\infty}$ is the symmetric positive definite matrix (\ref{def:Qinfty_rho_OU1}).
The operator $\mathfrak{L}=q^w(x,D_x)$ is a quadratic operator whose Weyl symbol
\begin{equation}\label{pav69}
q(x,\xi)=\frac{1}{2}|Q^{\frac{1}{2}}\xi|^2+\frac{1}{8}|Q^{\frac{1}{2}}Q_{\infty}^{-1}x|^2+i\Big\langle\Big(\frac{1}{2}QQ_{\infty}^{-1}+B\Big)x,\xi\Big\rangle, \quad (x,\xi) \in \rr^{2n},
\end{equation}
has a non-negative real part $\textrm{Re }q \geq 0$. This Weyl symbol is the complex conjugate of the Weyl symbol of the operator (\ref{pav6}). It follows that the singular space of the quadratic operator $\mathfrak{L}$ is zero $S=\{0\}$.  As in the previous section, the smallest integer $0 \leq k_0 \leq 2n-1$ satisfying 
\begin{equation}\label{pav309}
\Big(\bigcap_{j=0}^{k_0}\textrm{Ker}
\big[\textrm{Re }F(\textrm{Im }F)^j \big]\Big)\cap \rr^{2n}=\{0\},
\end{equation}
with $F$ being the Hamilton map of $q$,
corresponds exactly to the smallest integer $0 \leq k_0 \leq n-1$ satisfying
\begin{equation}\label{pav319}
\textrm{Rank}[Q^{\frac{1}{2}},BQ^{\frac{1}{2}},\dots, B^{k_0}Q^{\frac{1}{2}}]=n.
\end{equation}

Let $\omega$ be an open subset of $\mathbb{R}^n$ satisfying condition (\ref{hyp_omega}).
We can therefore deduce from Theorem~\ref{Main_result_control} applied to the quadratic operator $\mathfrak{L}$ that the parabolic equation 
\begin{equation}\label{evo2.91}
\left\lbrace\begin{array}{c}
\partial_th(t,x)+\mathfrak{L}h(t,x)=u(t,x)\un_{\omega}(x),\\
h|_{t=0}=h_0 \in L^2(\rr^n,dx),
\end{array}\right.
\end{equation}
is null-controllable from the set $\omega$ in any positive time $T>0$. Let $f_0 \in L_{1/\rho}^2$.
By using that the mappings 
\begin{equation}\label{evo2.9}
\begin{array}{cc}
\mathfrak{T} :  L^2 & \rightarrow  L^2_{1/\rho}\\
\ \ \ v &  \mapsto  \sqrt{\rho}v
\end{array}, \qquad \begin{array}{cc}
\mathfrak{T}^{-1} :  L^2_{1/\rho} & \rightarrow L^2\\
\quad \ \ \ v & \mapsto \sqrt{\rho}^{-1}v
\end{array},
\end{equation}
are isometric, we consider a control function $u \in L^2((0,T) \times \rr^n,dt \otimes dx)$ supported in $[0,T]\times\omega$ such that the mild solution $h$ to the equation (\ref{evo2.91}) with initial datum $h_0=\mathfrak{T}^{-1}f_0$ satisfies $h(T,\cdot)=0$. We deduce from (\ref{jen0772}), (\ref{pav4.59}) and (\ref{evo2.9}) that the mild solution
$$f=e^{-\frac{t}{2}\textrm{Tr}(B)}\mathfrak{T}h \in L_{1/\rho}^2,$$ 
to the equation
$$\left\lbrace \begin{array}{l}
\partial_t f(t,x) - \frac{1}{2}\textrm{Tr}[Q\nabla_x^2 f(t,x)]+\langle Bx, \nabla_x f(t,x)\rangle+\textrm{Tr}(B)f(t,x)=\tilde{u}(t,x)\un_{\omega}(x),\\
f|_{t=0}=f_0 \in L_{1/\rho}^2,                                       
\end{array}\right.$$
with the control function supported in $[0,T]\times\omega$,
$$\tilde{u}(t,x)=e^{-\frac{t}{2}\textrm{Tr}(B)}\mathfrak{T}u(t,x) \in L^2((0,T) \times \rr^n,dt \otimes \rho(x)^{-1}dx),$$
satisfies $f(T,\cdot)=0$. This proves that the hypoelliptic Fokker-Planck equation (\ref{syst_LR179}) is null-controllable from the set $\omega$ in any positive time $T>0$. This ends the proof of Corollary~\ref{cor_Main_result_control1}.

On the other hand, we deduce from Theorem~\ref{Main_result_obs} applied to the quadratic operator $\mathfrak{L}$ that for all $T>0$, there exists a positive constant $C_T>0$ such that 
\begin{equation}\label{bn1.08}
\forall h_0 \in L^2(\rr^n), \quad \| e^{-T\mathfrak{L}}h_0\|_{L^2(\rr^n,dx)}^2 \leq C_T\int_{0}^T \| e^{-t \mathfrak{L}}h_0\|_{L^2(\omega,dx)}^2dt.
\end{equation}
According to (\ref{pav4.59}), the semigroup $(e^{t\mathscr{P}})_{t \geq 0}$ on $L^2_{1/\rho}$ is given by 
\begin{equation}\label{bn1.8}
\forall f_0 \in L^2_{1/\rho}, \forall t \geq 0, \quad  e^{t\mathscr{P}}f_0=e^{-\frac{t}{2}\textrm{Tr}(B)}\mathfrak{T}e^{-t\mathfrak{L}}\mathfrak{T}^{-1}f_0,
\end{equation}
where $(e^{-t\mathfrak{L}})_{t \geq 0}$ denotes the $L^2(\rr^n,dx)$ contraction semigroup generated by $\mathfrak{L}$.  
By observing that 
$$\|e^{t\mathscr{P}}f_0\|_{L^2(\rr^n,\rho(x)^{-1} dx)}=e^{-\frac{t}{2}\textrm{Tr}(B)}\|e^{-t\mathfrak{L}}\mathfrak{T}^{-1}f_0\|_{L^2(\rr^n,dx)}$$
and 
$$\|e^{t\mathscr{P}}f_0\|_{L^2(\omega,\rho(x)^{-1} dx)}=e^{-\frac{t}{2}\textrm{Tr}(B)}\|e^{-t\mathfrak{L}}\mathfrak{T}^{-1}f_0\|_{L^2(\omega,dx)},$$
we deduce from (\ref{bn1.08}) and (\ref{bn1.8}) that the hypoelliptic Fokker-Planck operator $\mathscr{P}$ satisfies the observability estimate
$$\forall g \in L^2(\mathbb{R}^n,\rho(x)^{-1} dx), \quad \| e^{T\mathscr{P}} g\|_{L^2(\mathbb{R}^n,\rho(x)^{-1} dx)}^2 \leq \tilde{C}_T \int_0^T \| e^{t\mathscr{P}} g\|_{L^2(\omega,\rho(x)^{-1} dx)}^2 dt,$$
with $\tilde{C}_T=e^{-T\textrm{Tr}(B)}C_T>0$.
This ends the proof of Corollary~\ref{thm:OU1}.

\section{Application: Null-controllability and observability of a chain of two oscillators coupled to two heat baths at each side}\label{sec:appli}

This section is devoted to provide an application of the general results of null controllability and observability for accretive quadratic operators with zero singular spaces. 
This example given in~\cite[Section~4.3]{kps11} comes from the series of works~\cite{EH00,EH03,EPR99,HeHiSj2,HeHiSj3}. 
It is a model describing a chain of two oscillators coupled with two heat baths at each side. The particles are described by their respective positions and velocities $(x_j,y_j) \in \rr^{2d}$. For each oscillator, the particles are submitted to an external force derived from a real-valued potential $V_j(x_j)$ and a coupling between the two oscillators derived from a real-valued potential $V_c(x_2-x_1)$. We denote the full potential
$$V(x)=V_1(x_1)+V_2(x_2)+V_c(x_2-x_1), \quad x=(x_1,x_2) \in \rr^{2d},$$
$y=(y_1,y_2) \in \rr^{2d}$ the velocities and $z=(z_1,z_2) \in \rr^{2d}$ the variables describing the state of the particles in each of the heat baths. In each bath, the particles are submitted to a coupling with the nearest oscillator, a force given by the friction coefficient $\gamma$ and a thermal diffusion at the temperature $T_j$. We denote $w_1,w_2$ two standard $d$-dimensional Brownian motions and $w=(w_1,w_2)$. The system of equations describing this model is given by
\begin{equation}\label{osc1}
\left\lbrace
  \begin{array}{l}
    dx_1 = y_1 dt\\
    dx_2=y_2dt\\
    dy_1 = -\partial_{x_1}V(x)dt+z_1dt\\
    dy_2 = -\partial_{x_2}V(x)dt+z_2dt\\
    dz_1=- \gamma z_1 dt+\gamma x_1 dt-\sqrt{2 \gamma T_1}dw_1\\
    dz_2=-\gamma z_2dt+\gamma x_2 dt-\sqrt{2\gamma T_2}dw_2.
  \end{array}
\right.
\end{equation}
Setting $T_1=\frac{\alpha_1 h}{2}$ and $T_2=\frac{\alpha_2 h}{2}$, the corresponding equation for the density $g$ of particles is 
\begin{multline}\label{osc2}
h\partial_t g +\frac{\gamma}{2}\alpha_1 (-h\partial_{z_1})\Big(h\partial_{z_1}+\frac{2}{\alpha_1}(z_1-x_1)\Big)g \\ 
+ \frac{\gamma}{2}\alpha_2(-h\partial_{z_2})\Big(h\partial_{z_2}+\frac{2}{\alpha_2}(z_2-x_2)\Big)g + \Big(y \cdot h\partial_x-(\nabla_xV(x)-z)\cdot h\partial_y\Big)g =0.
\end{multline}
For simplicity, we consider the case when $h=1$, $\gamma=2$ and $d=1$. Furthermore, we consider the case when the external potentials are quadratic 
\begin{equation}\label{osc10}
V_1(x_1)=\frac{1}{2}ax_1^2, \quad V_2(x_2)=\frac{1}{2}b x_2^2, \quad V_c(x_1-x_2)=\frac{1}{2}c(x_1-x_2)^2,
\end{equation}
with $a,b,c \in \rr$. Let $\alpha>0$ be a positive constant. We assume that the parameters satisfy the following conditions
\begin{equation}\label{hyp:alpha_012}
\alpha > \frac{1}{2}\max(\alpha_1,\alpha_2), \quad \alpha_1>0,\quad \alpha_2>0, \quad (a+c-1)(b+c-1)-c^2 \neq 0.
\end{equation}
When these conditions hold, we can deduce from Theorems~\ref{Main_result_control} and~\ref{Main_result_obs} the following results:

\medskip

\begin{proposition}\label{propchain}
Let $\omega$ be an open subset of $\mathbb{R}_{x,y,z}^6$ satisfying (\ref{hyp_omega}). Setting
$$\Phi(x,y,z)=V(x)+\frac{|y|^2}{2}+\frac{|z|^2}{2}-z \cdot x, \qquad \mathcal{M}_{\alpha}=e^{-\frac{2\Phi}{\alpha }},$$
when the conditions (\ref{hyp:alpha_012}) hold, the evolution equation 
$$\left\lbrace
  \begin{array}{l}
\partial_t g(t,x,y,z)+\alpha_1 (-\partial_{z_1})\big(\partial_{z_1} +\frac{2}{\alpha_1}(z_1-x_1)\big)g(t,x,y,z) \\ 
\hspace{0.5cm} +\alpha_2(-\partial_{z_2})\big(\partial_{z_2}+\frac{2}{\alpha_2}(z_2-x_2)\big)g(t,x,y,z) + \big(y \cdot \partial_x-(\nabla_xV(x)-z)\cdot \partial_y\big)g(t,x,y,z)\\
 =u(t,x,y,z)\un_{\omega}(x,y,z),\\
g|_{t=0}=g_0 \in L^2(\mathbb{R}^6,\mathcal{M}_{\alpha}^{-1}dx dy dz),
 \end{array}
\right.$$
posed in the space $L^2(\mathbb{R}^6,\mathcal{M}_{\alpha}^{-1}dx dy dz)$ is null-controllable from the set $\omega$ in any positive time $T>0$, with a control function 
$u \in L^2((0,T) \times \rr^6, dt \otimes \mathcal{M}_{\alpha}^{-1}dx dy dz)$ supported in $[0,T]\times\omega$. On the other hand, the mild solution to the evolution equation 
$$\left\lbrace
  \begin{array}{l}
\partial_t g(t,x,y,z)+\alpha_1 (-\partial_{z_1})\big(\partial_{z_1} +\frac{2}{\alpha_1}(z_1-x_1)\big)g(t,x,y,z) \\ 
\hspace{0.5cm} +\alpha_2(-\partial_{z_2})\big(\partial_{z_2}+\frac{2}{\alpha_2}(z_2-x_2)\big)g(t,x,y,z) + \big(y \cdot \partial_x-(\nabla_xV(x)-z) \cdot \partial_y\big)g(t,x,y,z)=0,\\
g|_{t=0}=g_0 \in L^2(\mathbb{R}^6,\mathcal{M}_{\alpha}^{-1}dx dy dz),
 \end{array}
\right.$$
satisfies the following observability estimate: 
\begin{multline*}
\forall T>0, \exists C_T>0, \forall g_0 \in L^2(\mathbb{R}^6,\mathcal{M}_{\alpha}^{-1}dx dy dz), \\ \|g(T)\|_{L^2(\mathbb{R}^6,\mathcal{M}_{\alpha}^{-1}dx dy dz)}^2 \leq C_T \int\limits_0^T \|g(t)\|_{L^2(\omega,\mathcal{M}_{\alpha}^{-1}dx dy dz)}^2 dt.
\end{multline*}
\end{proposition}

\medskip

\begin{proof}
We consider the operator 
\begin{multline*}
P=\alpha_1 (-\partial_{z_1})\Big(\partial_{z_1} +\frac{2}{\alpha_1}(z_1-x_1)\Big) +\alpha_2(-\partial_{z_2})\Big(\partial_{z_2}+\frac{2}{\alpha_2}(z_2-x_2)\Big)\\ + \big(y \cdot \partial_x-(\nabla_xV(x)-z)\cdot \partial_y\big).
\end{multline*} 
We associate to the operator $P$ acting on $L^2(\mathbb{R}^6,\mathcal{M}_{\alpha}^{-1}dx dy dz)$ the quadratic operator $q^w(X,D_X)$ acting on $L^2(\rr^6,dxdydz)$ defined as
\begin{equation}\label{pav4.59ch}
\mathcal{M}_{\alpha}^{-\frac{1}{2}}P\mathcal{M}_{\alpha}^{\frac{1}{2}}=q^w(X,D_X)-2,
\end{equation}
with $X=(x,y,z) \in \rr^6$. The explicit computation of the quadratic operator
\begin{multline}\label{osc3}
q^w(X,D_X)=\alpha_1 \Big(-\partial_{z_1}+\frac{1}{\alpha}(z_1-x_1)\Big)\Big(\partial_{z_1}+\big(\frac{2}{\alpha_1}-\frac{1}{\alpha}\Big)(z_1-x_1)\Big) + \\ 
\alpha_2\Big(-\partial_{z_2}+\frac{1}{\alpha}(z_2-x_2)\Big)\Big(\partial_{z_2}+\Big(\frac{2}{\alpha_2}-\frac{1}{\alpha}\Big)(z_2-x_2)\Big) + 
\big(y\cdot\partial_x-(\nabla_xV(x)-z)\cdot\partial_y\big)+2,
\end{multline}
is led in~\cite[Section~4.3]{kps11}.
Its Weyl symbol is given by
\begin{multline*}
q=\alpha_1\zeta_1^2+\alpha_2\zeta_2^2+\beta_1(z_1-x_1)^2+\beta_2(z_2-x_2)^2+i\big[2\delta_1\zeta_1(z_1-x_1)+2\delta_2\zeta_2(z_2-x_2)\\ +y_1\xi_1+y_2\xi_2-\eta_1\big((a+c)x_1-cx_2-z_1\big)-\eta_2\big(-cx_1+(b+c)x_2-z_2\big)\big],
\end{multline*}
with 
$$\beta_1=\frac{\alpha_1}{\alpha}\Big(\frac{2}{\alpha_1}-\frac{1}{\alpha}\Big), \quad \beta_2=\frac{\alpha_2}{\alpha}\Big(\frac{2}{\alpha_2}-\frac{1}{\alpha}\Big), \quad \delta_1=\frac{\alpha_1}{\alpha}-1, \quad \delta_2=\frac{\alpha_2}{\alpha}-1,$$
where the notations $\xi, \eta,\zeta$ stand respectively for the dual variables of $x,y,z$. The condition 
$$\alpha > \frac{1}{2}\max(\alpha_1,\alpha_2),$$
ensures that this quadratic symbol has a non-negative real part $\textrm{Re }q \geq 0$. 
On the other hand, some algebraic computations led in~\cite[Section~4.3]{kps11} show that the Hamilton map $F$ of the quadratic symbol $q$ satisfies to
$$\textrm{Ker}(\textrm{Re }F) \cap \rr^{12}=\{(x,y,z,\xi,\eta,\zeta) \in \rr^{12} :\ \zeta=0, \ x=z\},$$
$$\textrm{Ker}(\textrm{Re }F) \cap \textrm{Ker}(\textrm{Re }F\textrm{Im }F) \cap \rr^{12}=\{(x,y,z,\xi,\eta,\zeta) \in \rr^{12} :\ y=\eta=\zeta=0, \ x=z\},$$
\begin{multline*}
\textrm{Ker}(\textrm{Re }F) \cap \textrm{Ker}(\textrm{Re }F\textrm{Im }F) \cap \textrm{Ker}\big(\textrm{Re }F(\textrm{Im }F)^2\big) \cap \rr^{12}\\ =\{ y=\xi=\eta=\zeta=0, \ x=z, \ (a+c-1)x_1-cx_2=0, -cx_1+(b+c-1)x_2=0\}.
\end{multline*}
According to (\ref{hyp:alpha_012}), the singular space is therefore equal to zero $S=\{0\}$.

Let $\omega$ be an open subset of $\mathbb{R}^6$ satisfying condition (\ref{hyp_omega}).
We can therefore deduce from Theorem~\ref{Main_result_control} applied to the quadratic operator $q^w(X,D_X)$ that the parabolic equation 
\begin{equation}\label{evo2.91ch}
\left\lbrace\begin{array}{c}
\partial_th(t,X)+q^w(X,D_X)h(t,X)=u(t,X)\un_{\omega}(X),\\
h|_{t=0}=h_0 \in L^2(\rr^6,dX),
\end{array}\right.
\end{equation}
is null-controllable from the set $\omega$ in any positive time $T>0$. Let $g_0 \in L^2(\mathbb{R}^6,\mathcal{M}_{\alpha}^{-1}dX)$.
By using that the mappings 
\begin{equation}\label{evo2.9ch}
\begin{array}{cc}
\mathfrak{T} :  L^2(\rr^6,dX)  & \rightarrow  L^2(\mathbb{R}^6,\mathcal{M}_{\alpha}^{-1}dX)\\
\ \ \ v  & \mapsto   \sqrt{\mathcal{M}_{\alpha}}v
\end{array}, \ \begin{array}{cc}
\mathfrak{T}^{-1} :  L^2(\mathbb{R}^6,\mathcal{M}_{\alpha}^{-1}dX) & \rightarrow L^2(\rr^6,dX)\\
\quad \ \ \ v & \mapsto \sqrt{\mathcal{M}_{\alpha}}^{-1}v
\end{array},
\end{equation}
are isometric, we consider a control function $u \in L^2((0,T) \times \rr^6,dt \otimes dX)$ supported in $[0,T]\times\omega$ such that the mild solution $h$ to the equation (\ref{evo2.91ch}) with initial datum $h_0=\mathfrak{T}^{-1}g_0$ satisfies $h(T,\cdot)=0$. We deduce from (\ref{pav4.59ch}) and (\ref{evo2.9ch}) that the mild solution
$$g=e^{2t}\mathfrak{T}h \in L^2(\mathbb{R}^6,\mathcal{M}_{\alpha}^{-1}dX),$$ 
to the equation
$$\left\lbrace
  \begin{array}{l}
\partial_t g(t,X)+\alpha_1 (-\partial_{z_1})\big(\partial_{z_1} +\frac{2}{\alpha_1}(z_1-x_1)\big)g(t,X) \\ 
\hspace{0.5cm} +\alpha_2(-\partial_{z_2})\big(\partial_{z_2}+\frac{2}{\alpha_2}(z_2-x_2)\big)g(t,X) + \big(y \cdot \partial_x-(\nabla_xV(x)-z)\cdot \partial_y\big)g(t,X)\\
 =\tilde{u}(t,X)\un_{\omega}(X),\\
g|_{t=0}=g_0 \in L^2(\mathbb{R}^6,\mathcal{M}_{\alpha}^{-1}dX),
 \end{array}
\right.$$
with the control function supported in $[0,T]\times\omega$,
$$\tilde{u}(t,X)=e^{2t}\mathfrak{T}u(t,X) \in L^2((0,T) \times \rr^6,dt \otimes \mathcal{M}_{\alpha}^{-1}dX),$$
satisfies $g(T,\cdot)=0$. It proves that this equation is null-controllable from the set $\omega$ in any positive time $T>0$.

On the other hand, we deduce from Theorem~\ref{Main_result_obs} applied to the quadratic operator $q^w(X,D_X)$ that for all $T>0$, there exists a positive constant $C_T>0$ such that 
\begin{equation}\label{bn1.08ch}
\forall h_0 \in L^2(\rr^6,dX), \ \| e^{-Tq^w(X,D_X)}h_0\|_{L^2(\rr^6,dX)}^2 \leq C_T\int_{0}^T \| e^{-tq^w(X,D_X)}h_0\|_{L^2(\omega,dX)}^2dt.
\end{equation}
According to (\ref{pav4.59ch}), the semigroup $(e^{-tP})_{t \geq 0}$ on $L^2(\mathbb{R}^6,\mathcal{M}_{\alpha}^{-1}dX)$ is given by 
\begin{equation}\label{bn1.8ch}
\forall g_0 \in L^2(\mathbb{R}^6,\mathcal{M}_{\alpha}^{-1}dX), \forall t \geq 0, \quad  e^{-tP}g_0=e^{2t}\mathfrak{T}e^{-tq^w(X,D_X)}\mathfrak{T}^{-1}g_0,
\end{equation}
where $(e^{-tq^w(X,D_X)})_{t \geq 0}$ denotes the $L^2(\rr^n,dx)$ contraction semigroup generated by $q^w(X,D_X)$.  
By observing that 
$$\|e^{-tP}g_0\|_{L^2(\mathbb{R}^6,\mathcal{M}_{\alpha}^{-1}dX)}=e^{2t}\|e^{-tq^w(X,D_X)}\mathfrak{T}^{-1}g_0\|_{L^2(\rr^6,dX)}$$
and 
$$\|e^{-tP}g_0\|_{L^2(\omega,\mathcal{M}_{\alpha}^{-1}dX)}=e^{2t}\|e^{-tq^w(X,D_X)}\mathfrak{T}^{-1}g_0\|_{L^2(\omega,dX)},$$
we deduce from (\ref{bn1.08ch}) and (\ref{bn1.8ch}) that the operator $P$ satisfies the observability estimate
$$\forall g \in L^2(\mathbb{R}^6,\mathcal{M}_{\alpha}^{-1}dX), \quad \| e^{-TP} g\|_{L^2(\mathbb{R}^6,\mathcal{M}_{\alpha}^{-1}dX)}^2 \leq \tilde{C}_T \int_0^T \| e^{-tP} g\|_{L^2(\omega,\mathcal{M}_{\alpha}^{-1}dX)}^2 dt,$$
with $\tilde{C}_T=e^{4T}C_T>0$.
This ends the proof of Proposition~\ref{propchain}.
\end{proof}

\section{Appendix}\label{appendix}

\subsection{Spectral inequality for Hermite functions}\label{refherm}

This section is devoted to the proof of a spectral inequality for Hermite functions. 
To that end, we use the following result proved by Le Rousseau and Moyano in~\cite[Proposition 3.2]{LeRousseau_Moyano}:

\medskip

\begin{proposition}\label{thm:poids}\emph{(Weight function for elliptic Carleman estimate).} 
Let $S > 0$, $Q= (0, S) \times \mathbb{R}^n$ and $\omega$ be an open subset of $\mathbb{R}^n$ satisfying (\ref{hyp_omega}).
There exists a function   $\psi \in  C^3([0,S]\times\mathbb{R}^n;\rr_+)$ such that
\begin{equation} \label{weight_0}
\psi \in W^{3,\infty}([0,S] \times \mathbb{R}^n)\,,
\end{equation}
\begin{equation} \label{weight_1}
\forall (s,x) \in Q, \quad |(\nabla_{s,x} \psi)(s,x)| \geq C, 
\end{equation}
\begin{equation} \label{weight_2}
\forall x \in \mathbb{R}^n \setminus \omega, \quad (\partial_s \psi) |_{s=0} \geq C,
\end{equation}
\begin{equation} \label{weight_3}
(\partial_s \psi) |_{s=S} \leq -C < 0\,, \quad \psi|_{s=S}=0\,,
\end{equation}
with $C > 0$ being a positive constant.
\end{proposition}

\medskip

In order to establish the spectral inequality for Hermite functions, we need to derive a global Carleman estimate for the augmented elliptic operator
$$\widetilde{P}= -\Delta_{s,x} + |x|^2 = -\partial_s^2-\Delta_x+|x|^2,$$
on the set $Q=(0,S)\times\mathbb{R}^n$.

\medskip

\begin{proposition} \label{thm:Carlm_ell}\emph{(Global elliptic Carleman estimate).} 
Let $\omega$ be an open subset of $\mathbb{R}^n$ satisfying (\ref{hyp_omega})
and $\psi$ the weight function given by Proposition~\ref{thm:poids}.
With
$$\varphi(s,x)=\exp(\lambda\psi(s,x)), \quad \lambda \geq 1,$$ 
there exist some positive constants $C>0$, $\tau_0\geq1$ and $\lambda_0 \geq 1$ such that
\begin{align}\label{estcarl}
& \ \ \tau^3 \|e^{\tau \varphi}g\|_{L^2(Q)}^2+\tau\|e^{\tau \varphi}xg\|_{L^2(Q)}^2+\tau \|e^{\tau \varphi}\nabla_{s,x}g\|_{L^2(Q)}^2\\ \notag
 +& \ \tau \|e^{\tau \varphi(0,\cdot)} (\partial_s g)|_{s=0}\|_{L^2(\mathbb{R}^n)}^2 +
\tau e^{2\tau}\|\partial_s g|_{s=S}\|_{L^2(\mathbb{R}^n)}^2+\tau^3 e^{2\tau}\| g|_{s=S}\|_{L^2(\mathbb{R}^n)}^2  \\ \notag
& \ \leq  C(\|e^{\tau \varphi} \widetilde{P}g\|_{L^2(Q)}^2 +\tau e^{2\tau}\|(\nabla_x g)|_{s=S}\|_{L^2(\mathbb{R}^n)}^2 \\ \notag
& \ +\tau e^{2\tau}\|(x g)|_{s=S}\|_{L^2(\mathbb{R}^n)}^2+\tau \|e^{\tau \varphi(0,\cdot)}(\partial_s g)|_{s=0}\|_{L^2(\omega)}^2),
\end{align}
with $\lambda=\lambda_0$, for all $\tau \geq \tau_0$ and $g \in C^2([0,S];\mathscr{S}(\mathbb{R}^n,\mathbb{C}))$ verifying $g|_{s=0}\equiv 0$.
\end{proposition}

\medskip

\begin{proof} 
The proof of this result is a slight adaptation of the proof of Proposition~3.3 given in~\cite{LeRousseau_Moyano}. 
Compared to the estimate appearing in~\cite{LeRousseau_Moyano} (Proposition~3.3), there are two additional terms in the estimate (\ref{estcarl}) coming from the quadratic potential $|x|^2$ in the operator $\widetilde{P}$. More specifically, there is the extra term 
$$\tau \| e^{\tau \varphi} x g\|_{L^2(Q)}^2,$$ 
appearing in the left-hand-side of the estimate (\ref{estcarl}) and the extra term 
$$\tau e^{2\tau}\|(x g)|_{s=S}\|_{L^2(\mathbb{R}^n)}^2,$$ 
appearing in its right-hand-side. In the following proof, we only emphasize the differences with the one given in~\cite[Proposition~3.3]{LeRousseau_Moyano}.

We observe that $\widetilde{P}=P+|x|^2$, where $P$ denotes the operator appearing in the proof given in~\cite[Proposition~3.3]{LeRousseau_Moyano}.
Keeping the very same notations as in~\cite{LeRousseau_Moyano}, the conjugated operator 
$$\widetilde{P}_\varphi=e^{\tau \varphi} \widetilde{P} e^{-\tau\varphi}, \quad \tau \geq 1,$$
can be written as 
$$\widetilde{P}_\varphi=A+i\widetilde{B},$$ 
with $A=A_1+A_2+A_3$, $\widetilde{B}=B_1+\widetilde{B}_2$, where the operators $A_1, A_2, B_1, \widetilde{B}_2$ are the same as the ones defined in~\cite[p. 3205]{LeRousseau_Moyano}, whereas the additional operator $A_3$ is given by $A_3=|x|^2$. Let $0<\mu <2$ be positive parameter. Following~\cite[p. 3206]{LeRousseau_Moyano}, we next write
$$\widetilde{P}_\varphi+\tau \mu \Delta \varphi=A+iB,$$
with $B=B_1+B_2$ and $B_2=-i\tau(1+\mu) \Delta \varphi$.
For $v \in C^2([0,S];\mathscr{S}(\mathbb{R},\mathbb{C}))$ verifying $v|_{s=0}=0$, by expanding the square of the norm
$$\|\widetilde{P}_\varphi v+\tau \mu \Delta \varphi \, v\|_{L^2(Q)}^2=\|Av+iBv\|_{L^2(Q)}^2,$$
we obtain from Proposition~\ref{thm:poids} that 
\begin{equation}\label{ssh1}
\textrm{Re}(Av,iBv)_{L^2(Q)}=\sum\limits_{\substack{1\leq j \leq 3\\ 1 \leq k \leq 2}}I_{j,k}
\lesssim \|\widetilde{P}_\varphi v \|_{L^2(Q)}^2 + \mathcal{O}_{\lambda}(\tau^2) \|v\|_{L^2(Q)}^2,
\end{equation}
where the terms $I_{j,k}$ are explicitly computed in \cite[formulas (3.13) to (3.16)]{LeRousseau_Moyano} for any $1 \leq j,k\leq 2$.
On the other hand, the new terms
$$I_{3,k}=\textrm{Re}(A_3v,iB_kv)_{L^2(Q)}, \qquad  1 \leq k \leq 2,$$
are given by
$$I_{3,2}=\tau(1+\mu)\int_{Q} |xv|^2 \Delta_{s,x} \varphi \, dx ds\,,$$
$$I_{3,1} = 2 \tau \textrm{Re}\Big(\int_Q |x|^2 v (\partial_s \varphi \ \partial_s \overline{v} + \nabla_x \varphi \cdot \nabla_x \overline{v}) dx ds \Big)= J_{3,1} + BT_{3,1},$$
with
\begin{multline*}
J_{3,1}=- \tau \int_{Q} |xv|^2 \partial_s^2 \varphi dx ds - \tau \int_Q |v|^2 \text{div}(|x|^2\nabla_x \varphi)dx ds\\
= - \tau \int_Q |xv|^2 \Delta_{s,x} \varphi\, dx ds - 2 \tau \int_Q |v|^2 x \cdot \nabla_x \varphi dx  ds,
\end{multline*}
\begin{equation}\label{ssh0}
BT_{3,1}=\tau \int_{\mathbb{R}^n} |x v|^2 \partial_s \varphi|_{s=S} dx,
\end{equation}
since $v|_{s=0}=0$.
Following~\cite{LeRousseau_Moyano}, we deduce from (\ref{ssh1}) that 
\begin{equation}\label{ssh2}
\textrm{Re}(Av,iBv)_{L^2(Q)}=\widetilde{J}+\widetilde{BT} \lesssim \|\widetilde{P}_\varphi v \|_{L^2(Q)}^2 + \mathcal{O}_{\lambda}(\tau^2) \|v\|_{L^2(Q)}^2,
\end{equation}
with $\widetilde{J}=J+J_{3,1}+I_{3,2}$, $\widetilde{BT}=BT+BT_{3,1}$, where the terms $J$ and $BT$ are defined in~\cite[p. 3207]{LeRousseau_Moyano}.
We first study the interior terms:

\medskip

\noindent \textit{Interior terms.} Following~\cite{LeRousseau_Moyano}, we have
\begin{equation}\label{ssh50}
\widetilde{J}=\int_Q ( \tau^3 \gamma_0 |v|^2 + \tau \gamma_1 |\nabla_{s,x} v|^2 + \tau \gamma_2 |xv|^2) dx ds + \widetilde{X}\,,
\end{equation}
with
$$\widetilde{X} = X - 2 \tau \int_Q |v|^2 x \cdot \nabla_x \varphi dx ds, \quad \gamma_2=\mu \Delta_{s,x} \varphi,$$
where the terms $\gamma_0$, $\gamma_1$ and $X$ are defined in~\cite[p. 3208]{LeRousseau_Moyano}. It is established in~\cite[formula (3.18)]{LeRousseau_Moyano} that 
\begin{equation}\label{ssh3}
\gamma_0 \gtrsim \lambda^4 \varphi^3, \quad  \gamma_1 \gtrsim \lambda^2 \varphi,
\end{equation}
when $\lambda \geq 1$ is sufficiently large.
On the other hand, we deduce from (\ref{weight_0}) and (\ref{weight_1}) that 
$$\gamma_2=\mu \Delta_{s,x} \varphi = \mu\lambda (\Delta_{s,x}\psi) \varphi + \mu\lambda^2 |\nabla_{s,x}\psi|^2 \varphi 
\gtrsim \lambda^2 \varphi,$$
when $\lambda \geq 1$ is sufficiently large.
It provides a new positive term in the left-hand side of the estimate of the type
\begin{equation}\label{ssh51}
\tau \lambda^2 \int_Q |xv|^2 \varphi dx ds \geq \tau \lambda^2 \int_Q |xv|^2 dx ds,
\end{equation}
since $\varphi \geq 1$, as $\psi \geq 0$.
We observe from Proposition~\ref{thm:poids}, (\ref{ssh50}), (\ref{ssh3}) and (\ref{ssh51}) that the additional term in $\widetilde{X}$ given by
$$- 2 \tau \int_Q |v|^2 x \cdot \nabla_x \varphi dx ds  = - 2 \tau \lambda \int_Q |v|^2 x \cdot (\nabla_x \psi)  \varphi dx ds,$$
can be absorbed 
$$\Big|2 \tau \int_Q |v|^2 x \cdot \nabla_x \varphi dx ds\Big| \lesssim \tau^{1/2}\lambda \int_Q|xv|^2 \varphi dxds+\tau^{3/2}\lambda \int_Q|v|^2 \varphi dxds,$$
when the parameters $\lambda$ and $\tau$ are sufficiently large by the following positive term 
$$ \int_Q ( \tau \lambda^2 |xv|^2 \varphi + \tau^3 \lambda^4 \varphi^3 |v|^2) dx ds,$$
since $\varphi \geq 1$, appearing in the estimate from below of the term $\widetilde{J}$.
By taking advantage of the estimate from below of the term $J$ in~\cite[formula (3.19)]{LeRousseau_Moyano} and by choosing the parameter $\lambda$ sufficiently large (fixed) and the parameter $\tau$ sufficiently large (arbitrary), we obtain that
\begin{equation}\label{ssh10}
\widetilde{J} \gtrsim \tau^3 \|v\|_{L^2(Q)}^2 + \tau  \|\nabla_{s,x} v \|_{L^2(Q)}^2 + \tau \|xv\|_{L^2(Q)}^2.
\end{equation}
We next study the boundary terms:

\medskip

\noindent 
\textit{Boundary terms.} It follows from Proposition~\ref{thm:poids} and (\ref{ssh0}) that 
$$BT_{3,1} =  \tau \lambda \int_{\mathbb{R}^n} |x v|^2 \partial_s \psi|_{s=S}  dx \gtrsim -\tau \lambda \int_{\mathbb{R}^n}|(x v)|_{s=S}|^2dx,$$
since $\varphi|_{s=S}=1$.
Putting together this estimate and the lower bound on the term $BT$ given in \cite[formula (3.23)]{LeRousseau_Moyano}, we obtain that there exist some positive constants $C_0,C_1>0$ such that for sufficiently large values of the parameters $\lambda \geq 1$ and $\tau \geq 1$,
\begin{multline}\label{ssh11}
\widetilde{BT} \geq  
 C_0(\tau^3 \lambda^3 \| v|_{s=S} \|_{L^2(\mathbb{R}^n)}^2  
+ \tau   \lambda   \| \partial_s v|_{s=S} \|_{L^2(\mathbb{R}^n)}^2 
+ \tau   \lambda   \| \varphi^{1/2} \partial_s v|_{s=0} \|_{L^2(\mathbb{R}^n)}^2) \\
 - C_1(\tau   \lambda \| \nabla_x v|_{s=S}\|_{L^2(\mathbb{R}^n)}^2 
+ \tau  \lambda \| xv|_{s=S} \|_{L^2(\mathbb{R}^n)}^2 
+ \tau   \lambda \| \varphi^{1/2} \partial_s v|_{s=0}\|_{L^2(\omega)}^2).  
\end{multline}
By collecting the estimates (\ref{ssh2}), (\ref{ssh10}) and (\ref{ssh11}) obtained for the interior and boundary terms, we deduce that 
$$\begin{array}{ll}
         & \tau^3 \|v\|_{L^2(Q)}^2 + \tau \|\nabla_{s,x} v \|_{L^2(Q)}^2 +\tau \|xv\|_{L^2(Q)}^2 + \tau^3 \| v|_{s=S} \|_{L^2(\mathbb{R}^n)}^2   \\
         & + \tau \| \partial_s v|_{s=S} \|_{L^2(\mathbb{R}^n)}^2 + \tau\| \partial_s v|_{s=0} \|_{L^2(\mathbb{R}^n)}^2 \\
\lesssim &
\|\widetilde{P}_\varphi v\|_{L^2(Q)}^2 + \tau^2 \|v\|_{L^2(Q)}^2 + 
\tau \|\nabla_x v|_{s=S}\|_{L^2(\mathbb{R}^n)}^2 + \tau\| xv|_{s=S} \|_{L^2(\mathbb{R}^n)}^2+\tau\|\partial_s v|_{s=0}\|_{L^2(\omega)}^2,
\end{array}$$
when the parameters $\lambda \geq 1$ (fixed) and $\tau \geq 1$ (arbitrary) are sufficiently large. For sufficiently large values of the parameter $\tau \geq 1$, it follows that 
$$\begin{array}{ll}
         & \tau^3 \|v\|_{L^2(Q)}^2 + \tau \|\nabla_{s,x} v \|_{L^2(Q)}^2 + \tau \|xv\|_{L^2(Q)}^2 + \tau^3 \| v|_{s=S} \|_{L^2(\mathbb{R}^n)}^2  \\
         & + \tau \|\partial_s v|_{s=S} \|_{L^2(\mathbb{R}^n)}^2 + \tau\|\partial_s v|_{s=0} \|_{L^2(\mathbb{R}^n)}^2 \\
\lesssim &
\|\widetilde{P}_\varphi v\|_{L^2(Q)}^2  + 
\tau \|\nabla_x v|_{s=S}\|_{L^2(\mathbb{R}^n)}^2 + \tau\| xv|_{s=S} \|_{L^2(\mathbb{R}^n)}^2 +  \tau \|\partial_s v|_{s=0}\|_{L^2(\omega)}^2.
\end{array}$$
We observe that all the above calculations still make sense when taking $v=e^{\tau \varphi} g$, with $g \in C^2([0,S];\mathscr{S}(\mathbb{R}^n,\mathbb{C}))$ verifying $g|_{s=0}\equiv 0$. By using classical arguments, we finally obtain the estimate (\ref{estcarl}).
This ends the proof of Proposition~\ref{thm:Carlm_ell}.
\end{proof}

We deduce from the global Carleman estimate derived in Proposition~\ref{thm:Carlm_ell} the proof of the spectral inequality for Hermite functions:

\medskip

\noindent \textbf{Proof of Proposition \ref{thm:Spect_Ineq}:}
Let $N \geq 1$ and $(b_\alpha)_{\alpha \in \mathbb{N}^n} \in \mathbb{C}^{\mathbb{N}^n}$.
We consider the function
\begin{equation}\label{estt100}
u(s,x)=\sum\limits_{|\alpha| \leq N} b_\alpha \psi_{\alpha}(x) \frac{\sinh(s \sqrt{2|\alpha|+n})}{\sqrt{2|\alpha|+n}}.
\end{equation}
This function belongs to the space $C^2([0,S];\mathscr{S}(\mathbb{R}^n,\mathbb{C}))$. It satisfies the conditions $u|_{s=0}\equiv 0$ and $\widetilde{P}u=(-\partial_s^2-\Delta_x+|x|^2)u=0$, 
since 
$$\forall \alpha \in \mathbb{N}^n, \quad (-\Delta_x+|x|^2)\psi_\alpha=(2|\alpha|+n)\psi_\alpha.$$
Applying the global Carleman estimate given in Proposition~\ref{thm:Carlm_ell} provides that for all $\tau \geqslant \tau_0$,
\begin{multline} \label{IS_1}
 \tau^2 \|u|_{s=S}\|_{L^2(\mathbb{R}^n)}^2 \\
\leq C \big(\|(\nabla_x u)|_{s=S}\|_{L^2(\mathbb{R}^n)}^2+\|(x u)|_{s=S}\|_{L^2(\mathbb{R}^n)}^2+e^{\tau M}\|(\partial_s u)|_{s=0}\|_{L^2(\omega)}^2\big),
\end{multline}
with $0 \leq M=2\exp(\lambda_0 \sup_{x \in \omega} \psi(0,x)) - 2 <+\infty$, where $\psi \in W^{3,\infty}([0,S] \times \mathbb{R}^n)$ is the non-negative weight function given in Proposition~\ref{thm:poids}. We observe that
\begin{equation}\label{estt101}
\left\| u|_{s=S}  \right\|_{L^2(\mathbb{R}^n)}^2 = \sum\limits_{|\alpha| \leq N} 
|b_\alpha|^2\Big(\frac{\sinh(S \sqrt{2|\alpha|+n} )}{\sqrt{2|\alpha|+n}}\Big)^2\,.
\end{equation}
On the other hand, by using the classical formula
$$\sqrt{2} \frac{\partial \psi_\alpha}{\partial x_j} = \sqrt{\alpha_j} \psi_{\alpha - e_j} - \sqrt{\alpha_j+1} \psi_{\alpha+e_j},$$
we deduce that for all $1 \leq j \leq n$,
$$\begin{array}{ll}
 & \|(\partial_{x_j} u)|_{s=S}\|_{L^2(\mathbb{R}^n)}^2 
= \Big\|\sum\limits_{|\alpha|\leq N} b_\alpha (\partial_{x_j} \psi_{\alpha})\frac{\sinh(S \sqrt{2|\alpha|+n} )}{\sqrt{2|\alpha|+n}}\Big\|_{L^2(\mathbb{R}^n)}^2\\
\leq &  \Big\|\sum\limits_{|\alpha|\leq N} b_\alpha  \sqrt{\alpha_j} \psi_{\alpha - e_j}\frac{\sinh(S \sqrt{2|\alpha|+n} )}{\sqrt{2|\alpha|+n}}\Big\|_{L^2(\mathbb{R}^n)}^2 + \Big\|\sum\limits_{|\alpha|\leq N} b_\alpha  \sqrt{\alpha_j+1} \psi_{\alpha+e_j}\frac{\sinh(S \sqrt{2|\alpha|+n} )}{\sqrt{2|\alpha|+n}}\Big\|_{L^2(\mathbb{R}^n)}^2\\
\leq & (2N+1) \sum\limits_{|\alpha|\leq N} |b_\alpha|^2  \Big(\frac{\sinh(S \sqrt{2|\alpha|+n} )}{\sqrt{2|\alpha|+n}}\Big)^2.
\end{array}$$
It follows that
\begin{equation} \label{IS2}
\|(\nabla_x u)|_{s=S}\|_{L^2(\mathbb{R}^n)}^2 \leq n(2N+1) \sum\limits_{|\alpha|\leq N} |b_\alpha|^2  \Big(\frac{\sinh(S \sqrt{2|\alpha|+n} )}{\sqrt{2|\alpha|+n}}\Big)^2.
\end{equation}
By using the other classical formula
$$\sqrt{2} x_j \psi_\alpha = \sqrt{\alpha_j+1} \psi_{\alpha+e_j} + \sqrt{\alpha_j} \psi_{\alpha-e_j},$$
we obtain by using the very same lines that
\begin{equation} \label{IS3}
\|(x u)|_{s=S}\|_{L^2(\mathbb{R}^n)}^2 \leq n(2N+1)\sum\limits_{|\alpha|\leq N} |b_\alpha|^2\Big(\frac{\sinh(S \sqrt{2|\alpha|+n} )}{\sqrt{2|\alpha|+n}}\Big)^2.
\end{equation}
We deduce from (\ref{estt100}), (\ref{IS_1}), (\ref{estt101}), (\ref{IS2}) and (\ref{IS3}) that for all $N \geq 1$, $\tau \geq \tau_0$, $(b_\alpha)_{\alpha \in \mathbb{N}^n} \in \mathbb{C}^{\mathbb{N}^n}$,
$$(\tau^2-2nC(2N+1)) \sum\limits_{|\alpha|\leq N} |b_\alpha|^2\Big(\frac{\sinh(S \sqrt{2|\alpha|+n})}{\sqrt{2|\alpha|+n}}\Big)^2 \leq
Ce^{\tau M} \Big\| \sum\limits_{|\alpha|\leq N} b_\alpha \psi_{\alpha}  \Big\|_{L^2(\omega)}^2\,.$$
By taking $\tau_N= \max\{ \tau_0 , \sqrt{2nC(2N+1)+1} \}$, we obtain that for all $N \geq 1$,
$$\sum\limits_{|\alpha|\leq N} |b_\alpha|^2 
\leq \frac{1}{S^2} \sum\limits_{|\alpha|\leq N} |b_\alpha|^2\Big(\frac{\sinh(S \sqrt{2|\alpha|+n} )}{\sqrt{2|\alpha|+n}}\Big)^2 
\leq \frac{Ce^{\tau_N M}}{S^2} \Big\| \sum\limits_{|\alpha|\leq N} b_\alpha \psi_{\alpha}  \Big\|_{L^2(\omega)}^2.$$
This ends the proof of Proposition~\ref{thm:Spect_Ineq}.

\subsection{Gelfand-Shilov regularity}\label{GSreg}

We refer the reader to the works~\cite{gelfand,rodino1,rodino,toft} and the references herein for extensive expositions of the Gelfand-Shilov regularity theory.
The Gelfand-Shilov spaces $S_{\nu}^{\mu}(\rr^n)$, with $\mu,\nu>0$, $\mu+\nu\geq 1$, are defined as the spaces of smooth functions $f \in C^{\infty}(\rr^n)$ satisfying the estimates
$$\exists A,C>0, \quad |\partial_x^{\alpha}f(x)| \leq C A^{|\alpha|}(\alpha !)^{\mu}e^{-\frac{1}{A}|x|^{1/\nu}}, \quad x \in \rr^n, \ \alpha \in \mathbb{N}^n,$$
or, equivalently
$$\exists A,C>0, \quad \sup_{x \in \rr^n}|x^{\beta}\partial_x^{\alpha}f(x)| \leq C A^{|\alpha|+|\beta|}(\alpha !)^{\mu}(\beta !)^{\nu}, \quad \alpha, \beta \in \mathbb{N}^n.$$
These Gelfand-Shilov spaces  $S_{\nu}^{\mu}(\rr^n)$ may also be characterized as the spaces of Schwartz functions $f \in \mathscr{S}(\rr^n)$ satisfying the estimates
$$\exists C>0, \eps>0, \quad |f(x)| \leq C e^{-\eps|x|^{1/\nu}}, \quad x \in \rr^n, \qquad |\widehat{f}(\xi)| \leq C e^{-\eps|\xi|^{1/\mu}}, \quad \xi \in \rr^n.$$
In particular, we notice that Hermite functions belong to the symmetric Gelfand-Shilov space  $S_{1/2}^{1/2}(\rr^n)$. More generally, the symmetric Gelfand-Shilov spaces $S_{\mu}^{\mu}(\rr^n)$, with $\mu \geq 1/2$, can be nicely characterized through the decomposition into the Hermite basis $(\Psi_{\alpha})_{\alpha \in \mathbb{N}^n}$, see e.g. \cite{toft} (Proposition~1.2),
\begin{multline*}
f \in S_{\mu}^{\mu}(\rr^n) \Leftrightarrow f \in L^2(\rr^n), \ \exists t_0>0, \ \big\|\big((f,\Psi_{\alpha})_{L^2}\exp({t_0|\alpha|^{\frac{1}{2\mu}})}\big)_{\alpha \in \mathbb{N}^n}\big\|_{l^2(\mathbb{N}^n)}<+\infty\\
\Leftrightarrow f \in L^2(\rr^n), \ \exists t_0>0, \ \|e^{t_0\mathcal{H}^{\frac{1}{2\mu}}}f\|_{L^2(\rr^n)}<+\infty,
\end{multline*}
where $\mathcal{H}=-\Delta_x+|x|^2$ stands for the harmonic oscillator.

\subsection{Adapted Lebeau-Robbiano method for observability}\label{Appendix_LM}

This appendix is written in collaboration with Luc Miller\footnote{Universit\'e Paris-Ouest, Nanterre La D\'efense, UFR SEGMI, B\^atiment G, 200 Av. de la R\'epublique, 92001 Nanterre Cedex, France (luc.miller@math.cnrs.fr)}, and provides a proof of Theorem~\ref{Meta_thm_AdaptedLRmethod}.

\begin{proof}
For simplicity, the notation $\|\cdot\|$ refers in all the following to the norm $\|\cdot\|_{L^2(\Omega)}$.

\medskip

\noindent
\emph{Step 1}: We begin by establishing the following estimate: $\forall q>0$, $\exists 0<\tau_0'(q)<t_0$, $\exists M(q)>0$,  
\begin{equation} \label{meta_thm_telescopic_formula}
\forall 0 < \tau <\tau_0'(q), \forall g \in L^2(\Omega), \quad f_q(\tau) \|e^{\tau A} g \|^2  - f_q(q\tau) \|g\|^2 \leq \int_{\frac{\tau}{2}}^\tau \|e^{tA} g \|_{L^2(\omega)}^2 dt,
\end{equation}
where 
\begin{equation}\label{meta_thm_def_f}
f_q(s)=\exp\Big(-\frac{M(q)}{s^{\frac{am}{b-a}}}\Big), \quad s>0.
\end{equation}
To that end, we consider
$$\forall q>0, \quad \gamma(q)= \Big( \frac{3 c_1 2^{a+m}}{c_2q^{\frac{am}{b-a}}}\Big)^{\frac{1}{b-a}}.$$
We observe that 
\begin{equation} \label{relation_gamma}
\forall q>0, \quad c_2 \gamma(q)^b 2^{-m} = 3 c_1 \big(2\gamma(q)\big)^a q^{-\frac{am}{b-a}}.
\end{equation}
For all $q>0$, we can find $0< \tau_0'(q) <t_0$ such that for all $0 < \tau < \tau_0'(q)$,
\begin{equation} \label{relation_tau0}
\frac{\gamma(q)}{\tau^{\frac{m}{b-a}}} > 1, \qquad
\frac{\tau}{4} \geq \exp\Big(- \frac{c_1 \big(2\gamma(q)\big)^a}{\tau^{\frac{am}{b-a}}} \Big), \qquad 
\frac{\tau}{c_2^2} \leq \exp\Big(\frac{c_2 \gamma(q)^b}{2^{m}\tau^{\frac{am}{b-a}}}\Big).
\end{equation}
Let $q>0$, $0< \tau < \tau_0'(q)$ and $g \in L^2(\Omega)$. 
There exists a positive integer $k(q,\tau) \geq 1$ verifying
\begin{equation} \label{relation_k}
1<\frac{\gamma(q)}{\tau^{\frac{m}{b-a}}} \leq k(q,\tau) \leq  \frac{2 \gamma(q)}{\tau^{\frac{m}{b-a}}},
\end{equation}
since according to (\ref{relation_tau0}), the interval $(\gamma(q) \tau^{-\frac{m}{b-a}}, 2 \gamma(q) \tau^{-\frac{m}{b-a}})$  is of length $>1$, and is contained in $(1,+\infty)$. 
We deduce from the Pythagorean identity, the triangular inequality and (\ref{Meta_thm_IS}) that for all $t>0$, $k \geq 1$ and $g \in L^2(\Omega)$,
\begin{multline*}
\frac{e^{- 2 c_1 k^a}}{2}\| e^{tA} g \|^2 
 \leq \frac{e^{- 2 c_1 k^a}}{2}( \| \pi_k e^{tA} g \|^2 + \| (1-\pi_k)e^{tA}g \|^2)  
 \leq \frac{1}{2}  \| \pi_k e^{tA} g \|_{L^2(\omega)}^2 +  \| (1-\pi_k) e^{tA} g \|^2  \\
 \leq \| e^{tA} g \|_{L^2(\omega)}^2 +  \| (1-\pi_k) e^{tA} g \|_{L^2(\omega)}^2 +  \| (1-\pi_k) e^{tA} g \|^2  \leq  \| e^{tA} g \|_{L^2(\omega)}^2 + 2 \| (1-\pi_k) e^{tA} g \|^2, 
\end{multline*}
since $\|\cdot\|_{L^2(\omega)} \leq \|\cdot\|$.
By integrating the previous estimate on the interval $(\frac{\tau}{2},\tau)$, and by using the contraction property of the semigroup, we deduce from (\ref{Meta_thm_IS}) and (\ref{Meta_thm_dissip}) that for all $k \geq 1$, $0<\tau <t_0$ and $g \in L^2(\Omega)$,
\begin{align} \label{meta_thm_interm1}
\frac{\tau}{4}\,  e^{- 2 c_1 k^a}\, \| e^{\tau A} g \|^2
& \leq \int_{\frac{\tau}{2}}^{\tau} \frac{1}{2}\, e^{- 2 c_1 k^a}\, \| e^{tA} g \|^2 dt \\ \notag
& \leq \int_{\frac{\tau}{2}}^{\tau} \big(\| e^{tA} g \|_{L^2(\omega)}^2 + 2 \| (1-\pi_k) e^{tA} g \|^2\big) dt \\ \notag
& \leq \int_{\frac{\tau}{2}}^{\tau}  \Big(\| e^{tA} g \|_{L^2(\omega)}^2 +  \frac{2}{c_2^2}\, e^{-2c_2t^m k^b }\, \| g \|^2\Big)dt \\ \notag
& \leq \int_{\frac{\tau}{2}}^{\tau}  \| e^{tA} g \|_{L^2(\omega)}^2 dt +  \frac{\tau}{c_2^2}\, e^{-2c_2(\frac{\tau}{2})^mk^b }\, \| g \|^2.
\end{align}
Setting $M(q)=3 c_1\big(2\gamma(q)\big)^a$, we deduce from (\ref{meta_thm_def_f}), (\ref{relation_gamma}), (\ref{relation_tau0}) and (\ref{relation_k}) that for all $q>0$, $0<\tau <\tau_0'(q)$,
\begin{multline} \label{meta_thm_interm2}
\frac{\tau}{4}\, e^{- 2 c_1 k(q,\tau)^a} \geq \frac{\tau}{4}\exp\Big(- \frac{2 c_1 (2\gamma(q))^a}{\tau^{\frac{am}{b-a}}}\Big)\\ \geq \exp\Big(- \frac{3 c_1 (2\gamma(q))^a}{\tau^{\frac{am}{b-a}} }\Big) = \exp\Big(- \frac{M(q)}{\tau^{\frac{am}{b-a}} }\Big) = f_q(\tau) 
\end{multline}
and
\begin{multline} \label{meta_thm_interm3}
\frac{\tau}{c_2^2}\, e^{-2c_2(\frac{\tau}{2})^mk(q,\tau)^b} \leq  \frac{\tau}{c_2^2}\exp\Big(- \frac{2 c_2 \gamma(q)^b}{2^{m}\tau^{\frac{am}{b-a}} }\Big) \leq \exp\Big(- \frac{c_2 \gamma(q)^b}{2^{m}\tau^{\frac{am}{b-a}}}\Big)\\
 = \exp\Big(-\frac{3c_1(2\gamma(q))^a}{(q\tau)^{\frac{am}{b-a}}}\Big) = \exp\Big(-\frac{M(q)}{(q\tau)^{\frac{am}{b-a}}}\Big)= f_q(q\tau).
\end{multline}
Then, the estimate (\ref{meta_thm_telescopic_formula}) readily follows from the estimates (\ref{meta_thm_interm1}), (\ref{meta_thm_interm2}) and (\ref{meta_thm_interm3}).

\bigskip

\noindent
\emph{Step 2}: We can now derive the observability estimate (\ref{meta_thm_IO}) from a telescopic serie argument 
due to \cite{miller2010} (see also~\cite{miller_SMF}) and already exploited in~\cite{apraiz,phung,YubiaoZhang}.

We consider the parameters $\tau_0'=\tau_0'(\frac{1}{2})$, $M=M(\frac{1}{2})$ and the function $f=f_{\frac{1}{2}}$ defined in Step~1 for the choice of parameter $q=\frac{1}{2}$.
We set
\begin{equation}\label{kl0}
C_1=M 2^{\frac{am}{b-a}}>0, \qquad \tilde{T}_0=2\tau_0'>0.
\end{equation}
For $0<T<\tilde{T}_0$, we define for all $k \geq 0$,
\begin{equation}\label{kl}
\tau_k=\frac{T}{2^{k+1}}\,, \qquad T_0=T\,, \qquad T_{k+1}=T_k-\tau_k.
\end{equation}
By applying the estimate (\ref{meta_thm_telescopic_formula}) to the function $e^{T_{k+1} A} g$ with the parameter $\tau_k$, we obtain that for all $k \geq 0$ and $g \in L^2(\Omega)$,
$$f(\tau_k) \|e^{T_k A} g \|^2 - f(\tau_{k+1}) \|e^{T_{k+1}A} g\|^2 \leq \int_{T_{k+1}}^{T_k} \| e^{tA} g \|_{L^2(\omega)}^2 dt\,.$$
Summing up the previous estimates for all $k \geq 0$ provides that
\begin{equation}\label{kl1}
f(\tau_0) \| e^{T A} g \|^2=f(\tau_0) \| e^{T_0 A} g \|^2  \leq \int_0^{T_0} \|e^{tA} g \|_{L^2(\omega)}^2 dt = \int_0^T \|e^{tA} g \|_{L^2(\omega)}^2 dt,
\end{equation}
since 
$$T_k \underset{k \rightarrow +\infty}{\longrightarrow} 0$$ 
and by the contractivity property of the semigroup
$$f(\tau_k) \|e^{T_k A} g \|^2 \leq \exp\Big(- \frac{M}{\tau_k^{\frac{am}{b-a}}}\Big) \|g\|^2 \underset{k \rightarrow +\infty}{\longrightarrow} 0.$$
We deduce from (\ref{kl0}), (\ref{kl}) and (\ref{kl1}) that 
\begin{equation}\label{kl3}
\forall 0<T<\tilde{T}_0, \forall g \in L^2(\Omega),  \quad  \| e^{T A} g \|^2 \leq \exp\Big(\frac{C_1}{T^{\frac{am}{b-a}}}\Big)\int_0^T \|e^{tA} g \|_{L^2(\omega)}^2 dt.
\end{equation}
Setting
$$C_2=\exp\Big(\frac{2^{\frac{am}{b-a}}C_1}{\tilde{T}_0^{\frac{am}{b-a}}}\Big)>1,$$
and by using anew the contractivity property of the semigroup, it follows from (\ref{kl3}) that for all $T \geq \tilde{T}_0$, $g \in L^2(\Omega)$,
\begin{multline}\label{kl4}
\| e^{T A} g \|^2 \leq \| e^{\frac{\tilde{T}_0}{2} A} g \|^2 \leq C_2\int_0^{\frac{\tilde{T}_0}{2}} \|e^{tA} g \|_{L^2(\omega)}^2 dt\\
 \leq C_2\exp\Big(\frac{C_1}{T^{\frac{am}{b-a}}}\Big)\int_0^{T} \|e^{tA} g \|_{L^2(\omega)}^2 dt.
\end{multline}
With $C=\sup(C_1,C_2)>1$, we deduce from (\ref{kl3}) and (\ref{kl4}) that 
\begin{equation}\label{kl5}
\forall T>0, \forall g \in L^2(\Omega),  \quad \| e^{T A} g \|^2 \leq  C\exp\Big(\frac{C}{T^{\frac{am}{b-a}}}\Big)\int_0^{T} \|e^{tA} g \|_{L^2(\omega)}^2 dt.
\end{equation}
This ends the proof of Theorem~\ref{Meta_thm_AdaptedLRmethod}.
\end{proof}

\end{document}